\documentclass[a4paper,10pt]{amsart}
\usepackage{amssymb,stmaryrd}
\usepackage{amsfonts}
\usepackage{amstext}
\usepackage{algorithmic}
\usepackage{algorithm}
\usepackage{young}
\usepackage{graphicx}
\usepackage{epstopdf}
\usepackage[all]{xy}
\usepackage{enumerate}
\usepackage{color}
\usepackage{epic}

\usepackage{tikz}
\usetikzlibrary{matrix,arrows,decorations.pathmorphing}

\numberwithin{equation}{section}
\newtheorem{theorem}{Theorem}[section]
\newtheorem{lemma}[theorem]{Lemma}
\newtheorem{proposition}[theorem]{Proposition}
\newtheorem{corollary}[theorem]{Corollary}

\theoremstyle{definition}
\newtheorem{definition}[theorem]{Definition}
\newtheorem{example}[theorem]{Example}

\theoremstyle{remark}
\newtheorem{remark}[theorem]{\bf{Remark}}

\newcommand{\R}{{\mathbb{R}}}

\newcommand{\C}{{\mathbb{C}}}

\newcommand{\Z}{{\mathbb{Z}}}

\newcommand{\A}{{\mathbb{A}}}

\newcommand{\<}{{\langle}}
\renewcommand{\>}{{\rangle}}

\newcommand{\CC}{{\mathcal{C}}}

\newcommand{\CF}{{\mathcal{F}}}

\newcommand{\CI}{{\mathcal{I}}}

\newcommand{\CL}{{\mathcal{L}}}
\newcommand{\CM}{{\mathcal{M}}}

\newcommand{\CR}{{\mathcal{R}}}

\newcommand{\CQ}{{\mathcal{Q}}}

\newcommand{\isom}{{\cong}}
\newcommand{\Ad}{{\rm Ad}}

\renewcommand{\ker}{{\rm{ker}}}

\newcommand{\tens}{\otimes}
\newcommand{\id}{{\rm id}}

\renewcommand{\o}{{}_{(1)}}
\renewcommand{\t}{{}_{(2)}}
\renewcommand{\th}{{}_{(3)}}
\newcommand{\fo}{{}_{(4)}}
\newcommand{\fiv}{{}_{(5)}}
\newcommand{\six}{{}_{\scriptscriptstyle(6)}}
\newcommand{\sev}{{}_{\scriptscriptstyle(7)}}
\newcommand{\ei}{{}_{\scriptscriptstyle(8)}}
\newcommand{\nine}{{}_{\scriptscriptstyle(9)}}
\newcommand{\ten}{{}_{\scriptscriptstyle(10)}}
\newcommand{\ele}{{}_{\scriptscriptstyle(11)}}
\newcommand{\twe}{{}_{\scriptscriptstyle(12)}}
\newcommand{\thir}{{}_{\scriptscriptstyle(13)}}

\newcommand{\rz}{{}_{(\bar 0)}}

\newcommand{\co}{{}_{(\bar 1)}}

\newcommand{\Bo}{{}_{\underline{(1)}}}
\newcommand{\Bt}{{}_{\underline{(2)}}}
\newcommand{\extd}{{\rm d}}
\newcommand{\del}{{\partial}}
\newcommand{\eps}{\epsilon}
\newcommand{\ev}{{\rm ev}}
\newcommand{\coev}{{\rm coev}}
\newcommand{\und}{\underline}

\newcommand{\la}{{\triangleright}}
\newcommand{\ra}{{\triangleleft}}

\newcommand{\lbiprod}{{>\!\!\!\triangleleft\kern-.33em\cdot}}
\newcommand{\rbiprod}{{\cdot\kern-.33em\triangleright\!\!\!<}}

\newcommand{\nquad}{{\!\!\!\!\!\!}}

\begin{document}

\title{Hodge Star as Braided Fourier Transform}

\author[S. Majid]{Shahn Majid}

\address{Queen Mary University of London\\
School of Mathematical Sciences, Mile End Rd, London E1 4NS, UK}
\email{s.majid@qmul.ac.uk}

\date{}

\begin{abstract} We study super-braided Hopf algebras $\Lambda$ primitively generated by finite-dimensional right crossed (or Drinfeld-Radford-Yetter) modules $\Lambda^1$ over a Hopf algebra $A$ which are quotients of the augmentation ideal $A^+$ under right multiplication and the adjoint coaction. Here super-bosonisation $\Omega=A\rbiprod\Lambda$ provides a bicovariant differential graded algebra on $A$. We introduce $\Lambda_{max}$ providing the maximal prolongation, while the canonical braided-exterior algebra $\Lambda_{min}=B_-(\Lambda^1)$ provides the Woronowicz exterior calculus. In this context we introduce a Hodge star operator $\sharp$  by super-braided Fourier transform on $B_-(\Lambda^1)$ and left and right interior products by braided partial derivatives. Our new approach to the Hodge star (a) differs from previous approaches in that it is canonically determined by the differential calculus  and (b) differs on key examples, having order 3 in middle degree on $k[S_3]$ with its 3D calculus and obeying the $q$-Hecke relation $\sharp^2=1+(q-q^{-1})\sharp$ in middle degree on $k_q[SL_2]$ with its 4D calculus. Our work also provided a Hodge map on quantum plane calculi and a new starting point for calculi on coquasitriangular Hopf algebras $A$ whereby any subcoalgebra $\CL\subseteq A$  defines a sub braided-Lie algebra and $\Lambda^1\subseteq \CL^*$ provides the required data $A^+\to \Lambda^1$. 
\end{abstract}

\maketitle

%\tableofcontents

\section{Introduction}

Differential exterior algebras $\Omega$ on quantum groups were extensively studied in the 1990s starting with \cite{Wor} and have a critical role as examples of noncommutative geometry more generally. However, one problem which has remained open since that era is the general construction of a Hodge star operator $\sharp$ in noncommutative geometry, even in the quantum group case. Until now the Hodge operator has been treated mainly in an ad-hoc manner in particular examples, motivated typically by $\sharp^2=\pm\id$ as a requirement, e.g. \cite{MaRai,GomMa:ele}. We also note a framework\cite{Heck} based on a pair of differential structures and contraction with a generalised metric, and \cite{Kus} in another q-deformation framework. By contrast our new approach depends canonically on the braided-Hopf algebra structure of the exterior algebra which applies at least for bicovariant calculi on quantum groups and covariant calculi on quantum-braided planes. Moreover, this new approach gives very different, and we think more interesting, answers than the previous approaches. Specifically, Section~3.2 includes the example of $k(S_3)$, the function algebra on the permutation group $S_3$, with its 2-cycle calculus, where $\sharp$ in middle degree obeys
\[ \sharp^3=\id\]
so that $\sharp$ as a whole has order 6. We also cover electromagnetism on this finite group as in \cite{MaRai} but using the new Hodge star and again achieving a full diagonalisation of the Laplace-Beltrami operator $\extd\delta+\delta\extd$.  Similarly, Section~4.2 compute our Hodge operator $\sharp$ for $k_q[SL_2]$ with its standard 4D calculus\cite{Wor} and in middle degree we find, rather unexpectedly, that it obeys the well-known $q$-Hecke relation
\[ \sharp^2=\id+(q-q^{-1})\sharp\]
when suitably normalised.  Both of these examples are very different from requiring $\sharp^2=\pm\id$. While the nicest version of the theory assumes a central bi-invariant metric and volume form, our Fourier approach is more general as illustrated in Section~3.1 on the quantum plane $\A_q^2$.

Conceptually, we adopt a novel point of view\cite{Ma:perm,Ma:self}  on what the Hodge star {\em is} even classically. Namely, at every point of a manifold $M$ of dimension $n$ the exterior algebra of differential forms has fibre the exterior algebra $\Lambda(\R^n)$ with generators $\theta_i=\extd x^i$ in local coordinates and the usual `Grassmann algebra' relations $\theta_i\theta_j+\theta_j\theta_i=0$. This is a finite-dimensional super-Hopf algebra and we can apply a super version of usual Hopf algebra Fourier transform $\CF:\Lambda\to \Lambda^*$ using the Berezin integral $\int {\rm Vol}=1$ where ${\rm Vol}=\theta_1\cdots\theta_n$. This then extends to the whole manifold and in the presence of a metric gives the classical Hodge operator $\sharp:\Omega^m\to \Omega^{n-m}$. This point of view was also recently used for the Hodge star on supermanifolds\cite{Cas}. The same approach applies to bicovariant differential exterior algebras on Hopf algebras, where we recall that these are parallelizable via an exterior algebra of left-invariant differential forms $\Lambda$ forming a super-braided-Hopf algebra, so we can do super-Fourier transform on this. 

In algebraic terms, $\Omega$ in the bicovariant case is a super-Hopf algebra equipped with a split super-Hopf algebra projection $\Omega\twoheadrightarrow H$. Hence by a super-version of the Radford-Majid theorem \cite{Rad, Ma:skl} one knows that $\Omega\isom A\rbiprod \Lambda$ is the super-bosonisation of $\Lambda$ as a super-braided Hopf algebra in the braided category of crossed (or Drinfeld-Radford-Yetter) modules \cite{Whi,Rad,Dri,Yet}, where we assume that $A$ has invertible antipode. Moreover, in the standard setting $\Omega$ is generated by degrees 0,1 and hence $\Lambda$ has primitive generators $\Lambda^1$. As a result we can focus on such primitively-generated super-braided Hopf algebras $\Lambda$ and translate elements of noncommutative geometry in these terms.  Particularly, the exterior derivative restricts to $\extd:\Lambda\to \Lambda$ to make $\Lambda$ into a differential graded algebra (but not within the category, $\extd$ is not a morphism) and in Section~2 we give an explicit construction of this in the case $\Lambda=\Lambda_{max}$ corresponding to the maximal prolongation of $\Omega^1$. Any other exterior algebra corresponds to a quotient of this and at the other extreme we revisit the more well-known case $\Lambda=\Lambda_{min}=B_-(\Lambda^1)$ given by the canonical (super) braided-linear space associated to an object of an Abelian braided category. This $B_\pm(\Lambda^1)$ construction appeared in the case where $\Lambda^1$ is rigid in \cite{Ma:dbos} while quadratic  primitively generated braided-Hopf algebras appeared in \cite{Ma:fre,Ma:book}. The $B_+$ construction as an algebra is often called a `Nichols algebra' cf\cite{Nic} or `Nichols-Woronowicz algebra' cf\cite{Wor} but we note that neither of these works considered $B_\pm$ as braided-Hopf algebras, that structure being introduced following the development of braided-Hopf algebras as above. The super-braided Hopf algebra interpretation of the Woronowicz construction of bicovariant differential exterior algebras was in \cite{Ma:perm,Ma:dcalc} among other works. In this context the universal property of $B_-$  corresponds in some sense to the minimal relations needed to ensure Poincar\'e duality, a remark that will be reflected in our approach to the Hodge star. It was also observed recently \cite{MaTao} that the exterior derivative $\extd$  on a bicovariant exterior algebra is not only a super-derivation but also a super-coderivation
\[ \Delta \extd=(\extd\tens \id + (-1)^{|\ |}\tens\extd)\Delta\]
where $|\ |$ is the degree operator and $\Delta$ is the super-coproduct. This turns out to be key to going the other way of building  $(\Omega,\extd)$ from data $\extd$ on $\Lambda$. Although these matters are somewhat familiar, the  braided approach to the differential structure requires  proofs which we provide as part of a necessary systematic treatment. Section~2.4 similarly provides a canonical construction for the differential exterior algebra on $B_+(\Lambda^1)$ as a quantum-braided plane. 

As part of this, and critical for Fourier transform on $B_\pm(\Lambda^1)$, is the notion of braided-exponential  in our approach to these algebras (being used notably in \cite{Ma:dbos} to inductivively build up the quasitriangular structure of quantum groups $U_q(g)$ as a succession of $q$-exponentials). Here \cite{Ma:dbos} 
\[ B_\pm(\Lambda^1)=T_\pm\Lambda^1/\oplus_m\ker[m;\pm\Psi]!\]
\[ \exp=\sum_m ([m, \pm\Psi]!^{-1}\tens\id)\coev_{\Lambda^1{}^{\tens m}}\]
are defined in terms of braided factorials $[m,\Psi]!$ as in \cite{Ma:fre,Ma:book}. We recall this theory in Section~2.2. In the $B_+$ case we have previously proposed braided Fourier theory on the Fomin-Kirillov algebra and its super-version as Hodge star in\cite{Ma:perm}, but without a systematic treatment.

Our central results appear in Section~3. If we have a unique bi-invariant top degree (say of degree $n$) then  super-braided-Hopf algebra Fourier transform gives us a map $\CF:\Lambda^m\to \Lambda^{*(n-m)}$ and in the presence of a quantum metric a Hodge star map $\sharp:\Lambda^m\to \Lambda^{n-m}$. This extends in the context above to the geometric $\sharp:\Omega^m\to \Omega^{n-m}$.  Proposition~\ref{outer} establishes in some generality that $\sharp$ commutes with the braided antipode $S$ and is involutive in degrees $0,1,n-1,n$. This in turn follows from some general results about super-braided Fourier transform in Section~3.1 which builds on our previous diagrammatic work, particularly\cite{KemMa}. Section~3.1  also covers the Hodge operator on the well-known quantum plane $\A_q^2$, where $\Lambda=B_-(\Lambda^1)=\A_q^{0|2}$, the fermionic quantum plane.

In Section~4 we focus on the case of coquastriangular Hopf algebras $(A,\CR)$ \cite{Dri,Ma:bg,Ma:book}. In line with the braided-Hopf algebra methods of the present paper, we first present a starting point for the construction of $\Lambda^1$ itself, namely as the dual of a braided-Lie algebra. We show that every subcoalgebra $\CL\subseteq A$ is a braided-Lie algebra in the sense introduced in \cite{Ma:blie}. This gives a significantly cleaner result than previous attempts such as \cite{GomMa} and builds on our recent work \cite{MaTao}. Everything is worked out in detail for $\C_q[SL_2]$ recovering previously `R-matrix' formulae when we take the standard matrix subcoalgebra, including the 4D braided-Lie algebra\cite{Ma:blie} and the Woronowicz 4D calculus\cite{Wor} from $\Lambda=\CL^*$. We then compute the canonical braided-Fourier transform on the latter with results as described above.

We work over a general ground field $k$ and $q\in k^\times$. In examples, we will assume characteristic zero for our calculations. We use the Sweedler notation $\Delta a=a\o\tens a\t$ for coproducts and $\Delta_R v=v\rz\tens v\co$ for right coactions (summations understood). We denote the kernel of the counit by $A^+$ and $\pi_\eps:A\to A^+$ defined by $\pi_\eps(a)=a-\eps(a)$ is the counit projection. We will make extensive use of the theory of braided-Hopf algebras \cite{Ma:bg} including the diagrammatic notation used for this kind of algebra in \cite{Ma:cro,Ma:bhop}.

\section{Braided construction of exterior algebras on Hopf algebras}

In this preliminary section we give a self-contained braided-Hopf algebra approach to bicovariant exterior algebras on Hopf algebras building on our  recent work \cite{MaTao}. We recall first that a differential graded algebra means a graded algebra $\Omega=\oplus_n\Omega^n$ equipped with a super-derivation $\extd$ increasing degree by 1 and squaring to 0. The standard setting is where $\Omega^1$ is spanned by elements of the form $a\extd b$ where $a,b\in A=\Omega^0$, and $\Omega$ is generated by degrees 0,1; in this case we say that we have an {\em exterior algebra} over $A$. Any first order differential structure $(\Omega^1,\extd)$ over $A$ can be extended to a maximal prolongation. 

\subsection{Maximal prolongation on a Hopf algebra} When $A$ is a Hopf algebra we can ask that left and right comultiplication extends to a bicomodule structure with coactions $\Delta_L,\Delta_R$ on $\Omega^1$ and $\extd$ a bicomodule map. In this case it is known that $\Omega$ becomes a super-Hopf algebra with coproduct that of $A$ on degree zero and $\Delta_L+\Delta_R$ on degree 1, see \cite{Brz, MaTao}. In this case $\pi:\Omega\to A$ which sends all degrees $>0$ to zero and is the identity on degree 0, forms a Hopf algebra projection split by the inclusion of $A$. As a result, assuming that the antipode of $A$ is invertible and by a super version of \cite{Rad,Ma:skl}, we have $\Omega\isom A\rbiprod \Lambda$ where $\Lambda$ is a super-braided Hopf algebra in the braided category of right $A$-crossed (or Drinfeld-Yetter) modules.  We recall that a right $A$-crossed module means a vector space $\Lambda^1$ which is both a right module and a right comodule such that
\[ \Delta_R(v\ra a)=v\rz\ra a\t\tens (Sa\o)v\co a\th,\quad\forall v\in \Lambda^1,\ a\in A.  \]
In this case there is an associated map $\Psi:\Lambda^1\tens \Lambda^1\to \Lambda^1\tens \Lambda^1$ defined by $\Psi(v\tens w)= w\rz\tens v\ra w\co$. A similar map for any pair of crossed modules makes the category of these braided when the antipode $S$ is invertible. Here $A^+=\ker\eps$ is itself a right crossed module by right multiplication and $\Ad_R(a)=a\t\tens (Sa\o)a\th$ and the result in  \cite{Wor} that first order calculi $(\Omega^1,\extd)$ are classified by ad-stable right-ideals can be viewed as saying that they are classified by surjective morphisms $\varpi:A^+\to \Lambda^1$. 

This point of view was recently used in \cite{MaTao} to generalise beyond the surjective case, where we do not assume that $\varpi$ is surjective. The exterior algebra is similarly given as bosonisation of a pair $(\Lambda,\extd)$ consisting of a primitively generated (by degree 1) super-braided Hopf algebra $\Lambda$ in the crossed module category equipped with a super-derivation (the restriction of $\extd$) which is a right $A$-comodule map and obeys $\extd^2=0$. This is required to have the further characteristic properties\cite{MaTao}
\begin{equation}\label{mc} \extd \varpi(a)+(\varpi\pi_\eps a\o)(\varpi\pi_\eps a\t)=0\end{equation}
\begin{equation}\label{dmod} (\extd \eta)\ra a-\extd(\eta\ra a)=(\varpi\pi_\eps a\o)\eta\ra a\t - (-1)^{|\eta|} (\eta\ra a\o) \varpi\pi_\eps a\t\end{equation}
\begin{equation}\label{dDelta} \und\Delta\extd\eta-(\extd\tens\id+(-1)^{|\ |}\tens\extd)\und\Delta\eta=(-1)^{|\eta\Bo|}\eta\Bo\rz\tens (\varpi\pi_\eps \eta\Bo\co)\eta\Bt\end{equation}
for all $a\in A^+$ and $\eta\in \Lambda$, where we underlined the braided-coproduct. It is shown in \cite{MaTao} that given such $(\Lambda,\extd)$, we obtain a bicovariant calculus $(\Omega,\extd)$ with $\extd a=a\o\varpi\pi_\eps a\t$ on degree 0 and that $\extd$ is also a supercoderivation. Here (\ref{mc}) is called the Maurer-Cartan equation cf\cite{Wor}. These results also clarify the surjective case:

\begin{lemma}\label{mclemma} In the case where $\varpi:A^+\to \Lambda^1$ is surjective, if $\extd$ is a super-derivation on $\Lambda$ with $\extd^2=0$ and obeys the Maurer Cartan equation then (\ref{dmod})-(\ref{dDelta}) hold. \end{lemma}
\proof First, it is straightforward (but in the 2nd case somewhat involved) to check that if (\ref{dmod})-(\ref{dDelta}) hold on $\omega$ then they hold on $\omega\eta$ for all $\eta\in\Lambda^1$. In both cases we use the super-derivation rule to break down $\extd(\omega\eta)$. For the second case we also need the super-braided-homomorphism property of the coproduct,  in the form
\[ \und\Delta(\omega\eta)=\omega\Bo\tens \omega\Bt\eta+(-1)^{|\omega\Bt|}\omega\Bo\eta\rz\tens \omega\Bt\ra\eta\co \]
We use ({\ref{dDelta}) on degree 1 in the form
\[ \und\Delta \extd \eta=\extd\eta\tens 1+1\tens\extd\eta-\eta\rz\tens \varpi\pi_\eps \eta\co\]
in the start of the induction and so as to be able to similarly compute $\und\Delta(\omega\extd\eta)$. We omit further details of the induction but we still need to establish both properties on degree 1. Thus
\begin{eqnarray*} (\extd \varpi(a))\ra b-\extd(\varpi(a)\ra b)&=&-((\varpi\pi_\eps a\o)(\varpi\pi_\eps a\t))\ra b+(\varpi\pi_\eps(a\o b\o ))(\varpi\pi_\eps(a\t b\t))\\
&=& (\varpi\pi_\eps b\o) \varpi(a b\t)+ (\varpi(a b\o))(\varpi \pi_\eps b\t)\\
&=& (\varpi\pi_\eps b\o) (\varpi (a)\ra  b\t)+ (\varpi(a)\ra b\o)(\varpi \pi_\eps b\t)
\end{eqnarray*}
for all $a\in A^+, b\in A$. We use $\pi_\eps(ab)=\pi_\eps(a)b+\eps(a)\pi_\eps b$ for all $a,b\in A$. For (\ref{dDelta}) we check the degree 1 version as
\begin{eqnarray*} \und\Delta \extd\varpi(a)&=&-\und\Delta((\varpi\pi_\eps a\o)(\varpi\pi_\eps a\t))\\
&=&\extd \varpi a\tens 1+1\tens\extd \varpi a-\varpi\pi_\eps a\o\tens \varpi\pi_\eps a\t-(\varpi\pi_\eps a\t)\rz\tens (\varpi\pi_\eps a\o)\ra(\varpi\pi_\eps a\t)\co\\
&=&\extd \varpi a\tens 1+1\tens\extd \varpi a-\varpi\pi_\eps a\o\tens \varpi\pi_\eps a\t-\varpi\pi_\eps a\th\tens \varpi((\pi_\eps a\o)Sa\t a\fo)\\
&=&\extd \varpi a\tens 1+1\tens\extd \varpi a- \varpi \pi_\eps a\t\tens \varpi\pi_\eps ((Sa\o)a\th) 
\end{eqnarray*}
for all $a\in A$. The latter part of the calculation here amounts to the identity
\begin{equation}\label{PsiDelta} \Psi(\varpi \pi_\eps\tens \varpi\pi_\eps)\Delta=(\varpi \pi_\eps\tens \varpi\pi_\eps)(\Delta-\Ad_R)\end{equation}
for the crossed module braiding. \endproof

Hence in order to construct $(\Lambda,\extd)$ in the surjective case it suffices to take (\ref{mc}) as a definition $\extd \varpi(a):=-(\varpi\pi_\eps a\o)(\varpi\pi_\eps a\t)$ and show that this is well-defined and extends as a super-derivation of square zero.  This will then make $(\Lambda,\extd)$ itself into a DGA over $k$. 

\begin{proposition}\label{maxpro} Let $\varpi: A^+\twoheadrightarrow \Lambda^1$ be a surjective morphism in the category of right crossed modules. Then
\[ \Lambda_{max}= T\Lambda^1/\< (\varpi\pi_\eps\tens \varpi\pi_\eps)\Delta\ker\varpi\>\]
together with $\extd$ defined by the Maurer-Cartan equation gives a super-braided Hopf algebra in the category which is also a differential graded algebra obeying (\ref{dmod})-(\ref{dDelta}).  Its bosonisation $\Omega_{max}=A\rbiprod \Lambda_{max}$ is the maximal prolongation differential calculus extending $(\Omega^1,\extd)$. 
\end{proposition}
\proof We quotient by the minimal subspace in degree 2 for which $\extd:\Lambda^1\to \Lambda^2$ is well-defined by the Maurer-Cartan equation. Let
\[ \del=\sum_{j=1}^{m}(-1)^{j+1}\Delta_j,\quad \del:A^{\tens m}\to A^{\tens(m+1)}\]
be the usual cobar coboundary, where $\Delta_j$ denotes the coproduct in the $i$'th position, and define $\extd$ on degree $m$ by $\extd .(\varpi\pi_\eps)^m=-\cdot (\varpi\pi_\eps)^{\tens (m+1)}\del$. This is well-defined for the same reason as before because $\del$ is the sum of terms each acting via the the coproduct. Clearly $\del$ is a super-derivation and squares to 0, so $\extd$ has the same features. We also need to check that we have a super-braided Hopf algebra. Since the algebra is quadratic, the main relation to check is
\[ \und\Delta ((\varpi\pi_\eps a\o)(\varpi\pi_\eps a\t)):=(\varpi\pi_\eps a\o)(\varpi\pi_\eps a\t)\tens 1+1\tens (\varpi\pi_\eps a\o)(\varpi\pi_\eps a\t)\]
\[\quad\quad\quad\quad\qquad\qquad+ \varpi\pi_\eps a\o\tens \varpi\pi_\eps a\t-\Psi(\varpi\pi_\eps a\o\tens \varpi\pi_\eps a\t)\]
vanlshes whenever $a\in \CI=\ker\varpi\pi_\eps$. This is clear for the first two terms and for the remaining two we use (\ref{PsiDelta}) to obtain $(\varpi\pi_\eps\tens\varpi\pi_\eps)\Ad_R(a)$ (much as in the proof of Lemma~\ref{mclemma}) which indeed vanishes as $\CI$ is Ad-stable because $\varpi$ was a morphism. Hence by the lemma we have the required data $(\Lambda_{max},\extd)$ and obtain a bicovariant calculus after bosonisation, something one can also check directly from $\CI$ an Ad-stable right ideal and the structure of $A\rbiprod \Lambda_{max}$.  It is also clear from the construction, since we imposed the minimal relations compatible with the Maurer-Cartan equation, that our calculus is isomorphic to the maximal prolongation.   \endproof

\subsection{Braided linear spaces}

Here we take an aside to recall the theory of braided-linear spaces introduced in \cite{Ma:fre,Ma:book} but in a cleaner form as braided operators rather than braided matrices. Braided linear spaces was our term for primitively generated graded braided Hopf algebras, with particular emphasis in \cite{Ma:dbos} on what have later been called  `Nichols-Woronowicz algebras'\cite{Baz}. If $V$ is an object in an Abelian braided category then it inherits  a morphism $\Psi=\Psi_{V,V}:V\tens V\to V\tens V$ obeying the braid relations. Our setting is categorical but we use only the pair $(V,\Psi)$ and tensor powers of $V$ in the following definition.

\begin{definition}\label{braidedbinom}\cite{Ma:fre,Ma:book} Let $(V,\Psi)$ be an object in an Abelian monoidal category and a braiding on it.  The braided binomials here are defined recursively by
\begin{gather*}
\left[ {n\atop r}; \Psi\right] =\Psi_r\Psi_{r+1}\cdots\Psi_{n-1}(\left[ {n-1\atop r-1};\Psi\right] \tens\id)+\left[ {n-1\atop r}; \Psi\right]\tens\id,\quad \left[ {n\atop 0}; \Psi\right]=\left[ {n\atop n}; \Psi\right]=\id.
\end{gather*}
where $0<r<n$ and $\Psi_i$ denotes $\Psi$ acting in the $i,i+1$ tensor factors. We also define `braided integers' 
\[ \left[n; \Psi\right]:=\left[ {n\atop 1}; \Psi\right]=\Psi_1\Psi_2\cdots\Psi_{n-1}+\left[ {n-1\atop 1}; \Psi\right]\tens\id=\id+\Psi_1+\Psi_1\Psi_2+\cdots+\Psi_1\Psi_2\cdots\Psi_{n-1}\] 
and `braided factorials' $[n,\Psi]!=(\id\tens [n-1,\Psi]!)[n,\Psi]$ where $[1,\Psi]!=\id_V$.  We take the convention $[0,\Psi]!=\id_{\und 1}$.
\end{definition}
These are operator versions of binomial coefficients and generalise $q$-binomials when applied to the category of $\Z$-graded vector spaces with braiding given by powers of $q$. Relevant to us, the braided factorials also generalise symmetrizers and antisymmetrizers. We need the following main theorem about them:
\begin{theorem}\label{binomfact}\cite{Ma:fre}\cite[10.4.12]{Ma:book}
\[ ([r;\Psi]!\tens[n-r;\Psi]!)\left[ {n\atop r}; \Psi\right] =[n,\Psi]!,\quad 0\le r\le n. \]
 \end{theorem}
 The proof in \cite{Ma:book} is written in matrix terms but immediately translates as operators in our setting. In fact, these results amount to identities in the group algebra of the braid group and are best done diagrammatically. The key observations are that 
 \begin{equation}\label{rr-1} ([r,\Psi]\tens\id)\Psi_r\cdots\Psi_{n-1}=\Psi_r\cdots\Psi_{n-1}([r-1,\Psi]\tens\id)+\Psi_1\cdots\Psi_{n-1}\end{equation}
since the first $r-1$ terms in $[r,\Psi]$ commute with $\Psi_r\cdots\Psi_{n-1}$ and 
\begin{equation}\label{functbinom} \Psi_1\cdots\Psi_{n-1}(\left[ {n-1\atop r}; \Psi\right] \tens\id)=(\id\tens\left[ {n-1\atop r}; \Psi\right] )\Psi_1\cdots\Psi_{n-1}\end{equation}
by functoriality (since the braided-binomial is a morphism) or directly by induction using Definition~\ref{braidedbinom} and repeated use of the braid relations.  Using these properties, \cite{Ma:fre, Ma:book} then proves  by induction on $n$ that
\[ ([r,\Psi]\tens\id)\left[ {n\atop r}; \Psi\right]=(\id\tens\left[ {n-1\atop r-1}; \Psi\right])[n,\Psi]\]
from which the theorem follows by repeated application. Also, by writing the above definitions as diagrams and turning the diagrams up-side down, we have co-binomial maps and co-integers defined by
\[ \left[ {n\atop r}; \Psi\right]' =(\id\tens\left[ {n-1\atop r-1};\Psi\right]')\Psi_1\cdots\Psi_{r-1}+\id\tens\left[ {n-1\atop r}; \Psi\right]',\quad \left[ {n\atop 0}; \Psi\right]'=\left[ {n\atop n}; \Psi\right]'=\id.\]
\[ [n,\Psi]'=1+\Psi_{n-1}+\Psi_{n-2}\Psi_{n-1}+\cdots+\Psi_1\cdots\Psi_{n-1}\]
Moreover,
\[  [n,\Psi]'!:=[n,\Psi]'([n-1,\Psi]'\tens\id)\cdots ([2,\Psi]'\tens\id)=[n,\Psi]!\]
\[ \left[ {n\atop r}; \Psi\right]'([r;\Psi]!\tens[n-r;\Psi]!) =[n,\Psi]!,\quad 0\le r\le n \]
where the factorials coincide by repeated use of the braid relations or because both cases can be written as $\sum_{\sigma\in S_n}\Psi_{i_1}\cdots\Psi_{i_{l(\sigma)}}$ where $\sigma=s_{i_1}\cdots s_{i_{l(\sigma)}}$ is a reduced expression in terms of simple transpositions $s_i=(i,i+1)$.  

Next we consider the tensor algebra $TV$ in an Abelian braided tensor category as a direct sum of different degrees $T_nV:=V^{\tens n}$ and product given by concatenation of $\tens$. Here $T_0V=k$ or more precisely the unit object of the category. The unit $\eta$ of the algebra $TV$ is the identity map from $k\to T_0V$. Thus
\[ (V\tens \cdots \tens V)\tens (V\tens \cdots\tens V)\to V\tens\cdots \tens V\]
is the identity map with suitable rebracketing (with $\Phi$ as necessary in the general case). We also consider the identity maps
\[ \eta_n: V^{\tens n}\to T_nV\]
with $\eta_0=\eta$. Although all these maps are the identity, we are viewing them in different ways. We will consider two different braided Hopf algebra structures $T_\pm V$ on $TV$, as a braided-Hopf algebra or as a super-braided Hopf algebra in the category. 

\begin{proposition}\label{TVcoprod}\cite{Ma:fre}\cite[Propn. 10.4.9]{Ma:book} The tensor algebra has a braided Hopf algebra/super Hopf algebra structure $T_\pm V$ with coproduct
\[\Delta_{T_nV}=\sum_{r=0}^n(\eta_r\tens\eta_{n-r})\circ\left[ {n\atop r}; \pm\Psi\right]. \]
for the two cases. The counit is $\eps_{T_n V}=0$ for all $n>0$. 
 \end{proposition} 
 \proof This is the content of \cite[Ma:fre,Propn 10.4.9]{Ma:book} in the free case where we impose no relations, but we rework the proof in the current more formal notations and we state the super case explicitly. We start with  the linear coproduct
\[ \Delta_{T_1V}=\id_V\tens \eta+\eta\tens \id_V\]
and for $T_+V$ we extend this as a Hopf algebra in the braided category, while for $T_-V$ we extend as a super-Hopf algebra in the braided category. We do the first case; the other is exactly the same by replacing $\Psi$ by $-\Psi$.  The proof is by induction assuming the formula for $\Delta_{T_{n-1}V}$,
\begin{eqnarray*} \Delta_{T_nV}&=&(\Delta_{T_{n-1}V}).(\id_V\tens\eta+\eta\tens \id_V)=\left(\eta\tens\id_{V^{\tens n-1}}+\sum_{r=1}^{n-1}\left[ {n\atop r}; \Psi\right]\right)\und\cdot(\id_V\tens\eta+\eta\tens \id_V)\\
&=&\eta\tens\id_{V^{\tens n}}+(\eta_1\tens\eta_{n-1})\Psi_1\cdots\Psi_{n-1}+\sum_{r=1}^{n-1}(\eta_r\tens\eta_{n-r})(\left[ {n-1\atop r};\Psi\right]\tens\id_V)\\
&&+ \sum_{r=1}^{n-1}(\eta_{r+1}\tens\eta_{n-1-r})\Psi_{r+1}\cdots\Psi_{n-1}\left[ {n-1\atop r}; \Psi\right]+\id_{V^{\tens n}}\tens\eta\\
&=&\eta\tens\id_{V^{\tens n}}+\id_{V^{\tens n}}\tens\eta\\
&&+\sum_{r=1}^{n-1}(\eta_r\tens\eta_{n-r})(\left[ {n-1\atop r}; \Psi\right]\tens\id_V)+ \sum_{r=1}^{n-1}(\eta_{r}\tens\eta_{n-r})\Psi_{r}\cdots\Psi_{n-1}\left[ {n-1\atop r-1}; \Psi\right]\\
&=&\eta\tens\id_{V^{\tens n}}+\id_{V^{\tens n}}\tens\eta+\sum_{r=1}^{n-1}(\eta_r\tens\eta_{n-r})\left[ {n\atop r}; \Psi\right]=\sum_{r=0}^n(\eta_r\tens\eta_{n-r})\circ\left[ {n\atop r}; \Psi\right]\\
\end{eqnarray*}
where in the first line we split off the $r=0$ part of $\Delta_{T_{n-1}V}$ and $\und\cdot$ is the braided tensor product. We then compute out the latter and for the 4th equality we renumber $r+1\mapsto r$ in the second sum and absorb the otherwise missing $r=1$. For the 5th equality we use Definition~\ref{braidedbinom} and finally combine terms to obtain the desired expression for $\Delta_{T_n}V$.  Again, this is really a result at the level of the braid group algebra and can be done with diagrams. \endproof
 
In this situation we are now ready to define the (super)Hopf algebra quotients 
\begin{equation}\label{BV} B_\pm(V)=T_\pm V/\oplus_m \ker [m,\pm\Psi]!\end{equation}
as the {\em braided-symmetric algebra} and {\em braided exterior algebra} on $V$ respectively. That the coproduct descends to $B_\pm(V)$ follows immediately from Proposition~\ref{TVcoprod} and Theorem~\ref{binomfact}. That $\oplus_m[m,\pm\Psi]!$ is an ideal or equivalently that the product in $TV$ descends to the quotient follows from the arrow-reversed version of Theorem~\ref{binomfact} where the factorials are on the right.  It is easy to see that when our construction is in a braided category and $\phi:V\to W$ is a morphism then $\phi^{\tens}$ (the relevant power) in each degree is a morphism $B_\pm(V)\to B_\pm(W)$ of (super)braided-Hopf algebras. This is because, by functoriality of the braidings, the braidings and braided factorials are intertwined by $\phi$ on each strand in the diagrammatic picture. The $B_+(V)$ case  is also called the Nichols-Woronowicz algebra of $V$ due to the structure of the algebra, but the above description and the fact that it is a (super) braided-Hopf algebra is due to the author. The earliest examples were the braided-line and braided quantum-plane (see \cite{Ma:fre,Ma:book}) while other early examples were $U_q(n_+)$ in the work of Lusztig\cite{Lus}.  

Our own motivation to consider (\ref{BV}) to all degrees of relations was in the case  when $V$ has a right dual $V^*$. Recall that $V^{*\tens n}$ is right-dual to $V^{\tens n}$ by the nested use of $\ev_V$ and we use the same nesting convention for a duality pairing $\<\ ,\ \>$ on tensor products. In this case \cite{Ma:fre,Ma:book}  the tensor algebras $T_\pm V^*$  and $T_\pm V$ are dually paired by 
\begin{equation}\label{TVeval} \<\ , \>|_{T_nV^*\tens T_m V}=\delta_{n,m}\ev_{V^{\tens n}}(\id\tens [n,\pm\Psi]!) \end{equation}
and $B_\pm(V^*),B_\pm(V)$ are  clearly the quotients by the kernel of the pairing. 
That the product on one side is the coproduct on the other follows immediately from Proposition~\ref{TVcoprod} and Theorem~\ref{binomfact}. This means that $B_\pm(V^*), B_\pm(V)$ are nondegenerately paired (super) Hopf algebras in the braided category and the relations of $B_\pm(V)$ are the minimal relations compatible with this duality. The merit of this approach is that we also have the immediate result, which will need later:

\begin{corollary}\label{coev}\cite{Ma:dbos} If $V$ has a right dual and $B_\pm(V)$ has a finite top degree then it has a right dual via $\ev=\<\ ,\ \>$ and  coevaluation map $\coev:\und 1\to B_\pm(V)\tens B_\pm(V^*)$ is given by
\[ \exp_V:=\coev=\sum_m ( [m,\pm\Psi]!^{-1}\tens\id)\coev_{V^{\tens m}}\]
where the construction is independent of the choice of inverse image of $[m,\pm\Psi]!$. 
\end{corollary}

This makes more precise the notion of braided-exponentials in \cite{Ma:fre,Ma:book} without formally assuming that the braided factorials are invertible. It was used explicitly in \cite{Ma:dbos} to construct quasitriangular structures.  If we take the well-known case $\Lambda^1=kx$ in the braided category of $\Z/(n+1)$-graded vector spaces with braiding $\Psi(v\tens w)=q^{\deg(w)\deg(v)}w\tens v$, where $x$ has degree 1, and $q$ a primitive $n+1$-th root of $1$, we have $B_+(kx)=k[x]/(x^{n+1})$ and $\exp$ here is the truncated $q$-exponential where $m!$ is replaced by $[m,q]!$ and $[m,q]=(1-q^m)/(1-q)$ are $q$-integers.

\subsection{Minimal prolongation on a Hopf algebra}

We now return to our setting of a Hopf algebra $A$ with invertible antipode and a surjective morphism $\varpi:A^+\to \Lambda^1$ in the braided category of right $A$-crossed modules. The following braided-Hopf algebra version of Woronowicz's construction\cite{Wor} is largely known eg \cite{Ma:dcalc} but we provide a new direct construction for $\extd$ on $B_-(\Lambda^1)$ going  through (\ref{mc})-(\ref{dDelta}) from \cite{MaTao}. This is a rather different from the approach in \cite{Wor}, which was to formally adjoin an inner element $\theta$ and define $\extd=[\theta, \ \}$. 

\begin{proposition}\label{bicovextalgthm}  Let $\varpi:A^+\to \Lambda^1$ be a surjective morphism is the category of right crossed modules and $\Lambda_{min}=B_-(\Lambda^1)$. This is a quotient of $\Lambda_{max}$ and inherits $\extd$ obeying (\ref{mc})-(\ref{dDelta}). Its super-bosonisation $\Omega_{min}=A\rbiprod \Lambda_{min}$  recovers the Woronowicz bicovariant calculus\cite{Wor} on $A$ associated to the Ad-stable right ideal $\ker\varpi$. \end{proposition}
\proof We start with the identity (\ref{PsiDelta}) and Ad-invariance of $\CI$ implies now that $(\varpi\pi_\eps \tens \varpi\pi_\eps)\Delta(\CI)\subset \ker[2,-\Psi]=\id-\Psi$, meaning that the relations $(\varpi\pi_\eps \tens \varpi\pi_\eps)\Delta(\CI)=0$ of the maximal prolongation in Proposition~\ref{maxpro} already hold among the quadratic relations in $\Lambda_{min}=B_-(\Lambda^1)$. The latter is therefore a quotient of the $\Lambda_{max}$. Next we consider the coboundary $\del_m: A^{\tens m}\to A^{\tens(m+1)}$ as the proof of Proposition~\ref{maxpro} and regard $A$ as a right crossed module by $\Ad_R$ and $a\ra b=\pi_\eps(a)b$. Then $\pi_\eps$ becomes a morphism.  Lemma~\ref{PsiDeltalemma} below shows that $\del$ descends to $B_-(A)\to B_-(A)$ in each degree as it respects the kernels of the relevant braided-factorials. Next, $\pi_\eps$ being a crossed module morphism induces by $\pi_\eps^{\tens m}$ in degree $m$ a map $B_-(A)\to B_-(A^+)$, under which $\del$ descends to a map $-\tilde \extd:B_-(A^+)\to B_-(A^+)$ given by $-\tilde \extd\pi_\eps^{\tens m}=\pi_\eps^{\tens m}\del_m$ because the kernel of $\pi_\eps^{\tens m}$ is spanned by elements where at least one of the tensor factors is 1. When we apply $\del_m$ then every term has at least one tensor factor $1$ which is then killed by the final $\pi_\eps^{\tens m}$. This is the first cell of 
\[
\xymatrix{
B_-(A)  \ar[d]^{\del} \ar[r]^{\pi_\eps^{\tens}}  &  B_-(A^+)  \ar[d]^{-\tilde \extd} \ar[r]^{\varpi^{\tens}}   & B_-(\Lambda^1) \ar[d]^{-\extd}    \\
B_-(A)  \ar[r]^{\pi_\eps^{\tens}}   &  B_-(A^+) \ar[r]^{\varpi^{\tens}}     &    B_-(\Lambda^1) 
  } 
\]
Similarly, $\varpi:A^+\to \Lambda^1$ being a morphism of crossed modules induces $B_-(A^+)\to B_-(\Lambda^1)$ given by $\varpi^{\tens n}$ in degree $n$, and $\tilde\extd$ descends to a map $\extd:B_-(\Lambda^1)\to B_-(\Lambda^1)$ defined by $\extd \varpi^{\tens m}=\varpi^{\tens m}\tilde\extd$. This is because the kernel of $\varpi^{\tens m}$ consists of terms where at least one of the tensor factors is in $\CI$. When we compute the $\tilde\extd$ of such terms using $\del_m$, either a $\Delta_j$ does not act on this tensor factor, in which case this tensor factor is present in the output of $\Delta_j$ and the whole term is killed by the action of $(\varpi\pi_\eps)^{\tens m}$, or $\Delta_j$ does act on this element. But then $\cdots\tens (\varpi\pi_\eps\tens\varpi)\pi_\eps(\Delta \CI)\tens\cdots$ is in the kernel of $\id-\Psi$ in the relevant place as seen above, hence vanishes in $B_-(\Lambda)$. We are using the fact that the kernel in each degree contains the degree 2 relations between adjacent tensor factors. In this way, $\extd$ equips $\Lambda_{min}=B_-(\Lambda^1)$ with a differential as a quotient of the construction for $\Lambda_{max}$ in Proposition~\ref{maxpro}.  One can show that if $\Omega^1$ is inner by $\theta\in \Lambda^1$ then the same applies to $\Omega$, but we are not assuming this. The algebra structure of the bosonosation $\Omega=A\rbiprod\Lambda_{min}$ is more well-known to be isomorphic to the one in \cite{Wor}.  \endproof
 
 The following lemma was needed to complete the proof. Here $\Psi$ is the braiding for the crossed-module structure on $A$ whereby $\pi_\eps$ becomes a morphism. 
 
 \begin{lemma}\label{PsiDeltalemma} Let $A$ be a Hopf algebra and $\Psi_i=\Psi$ the induced braiding 
 \[ \Psi(a\tens b)=b\t\tens a (Sb\o)b\th-\eps(a)b\t\tens (Sb\o)b\th\]
acting in the $i,i+1$ position of a tensor power, $\Ad_i=\Ad_R$, $\Delta_i=\Delta$ in the $i$'th position. Then
 \[ [m,-\Psi]!\del_{m-1}=\left(\sum_{j=1}^{m-1}(-1)^{j+1}\Ad_1^{\tens j}\right)[m-1,-\Psi]!\]
 Here $\Ad^{\tens j}$ denotes the tensor product right coaction on $j$ copies (acting here in the first position). 
\end{lemma}
\proof Clearly
\[ \Psi_i\Delta_j=\Delta_j\Psi_i,\quad {\rm if}\  i<j-1,\quad  \Psi_i\Delta_j=\Delta_j\Psi_{i-1},\quad {\rm if}\ i>j+1\]
since the operators act on different tensor factors, just the numbering changes in the 2nd case. We also find by direct computation in the Hopf algebra that 
\[ \Psi_i\Delta_i=\Delta_i-\Ad_i,\quad  \Psi_i\Delta_{i-1}= (\Delta\tens\Ad)_{i-1}-{\Ad}_i,\quad \Psi_i(\Delta\tens\Ad)_i=(\Delta_{i+1}-\Ad_i)\Psi_i\]
where $\Ad=\Ad_R$ and $\Delta\tens\Ad$ is the tensor product right coaction. As a warm-up, using these relations,  we show
\begin{eqnarray*} [3,-\Psi]\del_2&=&(1-\Psi_1+\Psi_1\Psi_2)(\Delta_1-\Delta_2)\\
&=&\Delta_1-\Delta_2-\Psi_1\Delta_1+\Psi_1\Delta_2-\Psi_1\Psi_2(\Delta_2-\Psi_1)\\
&=&\Delta_1-\Delta_2-\Delta_1+\Ad_1+\Psi_1\Delta_2-\Psi_1(\Delta_2-(\Delta\tens\Ad)_1)\\
&=&-\Delta_2+\Ad_1+(\Delta_2-\Ad_1)\Psi_1=(\Ad_1-\Delta_1+\del_2)[2,-\Psi]
\end{eqnarray*}
Starting with this, we next prove by induction that
\begin{equation}\label{PsiDeltarelns}[m,-\Psi]\del_{m-1}=(\Ad_1-\Delta_1+\del_{m-1})[m-1,-\Psi].\end{equation}
Assuming this for $m-1$ in the role of $m$,  for the 2nd equality, 
\begin{eqnarray*} [m,-\Psi]\del_{m-1}&=&[m-1,-\Psi]\del_{m-2}+[m-1,-\Psi](-1)^m\Delta_{m-1}+(-1)^{m-1}\Psi_1\cdots\Psi_{m-1}\del_{m-1}\\
&=&(\Ad_1-\Delta_1+\del_{m-2})[m-2,-\Psi]+(-1)^m\Delta_{m-1}[m-2,-\Psi]+\Psi_1\cdots\Psi_{m-2}\Delta_{m-1}\\
&&-\Psi_1\cdots\Psi_{m-1}\Delta_{m-1}+ (-1)^{m-1}\Psi_1\cdots\Psi_{m-1}\del_{m-2}\\
&=&(\Ad_1-\Delta_1+\del_{m-1})[m-2,-\Psi]+\Psi_1\cdots\Psi_{m-2}\Ad_{m-1}+ (-1)^{m-1}\Psi_1\cdots\Psi_{m-1}\del_{m-1} \end{eqnarray*}
where we picked out and computed the $\Psi_1\cdots\Psi_{m-1}\Delta_{m-1}$ term from the sum in $\del_{m-1}$. Looking now at the 
last expression, we compute
\begin{eqnarray*}\Psi_1\cdots\Psi_{m-1}\del_{m-2}&=&\sum_{j=1}^{m-2}(-1)^{j+1}\Psi_1\cdots\Psi_{j+1}\Delta_j\Psi_{j+1}\cdots\Psi_{m-2}\\
&=&\sum_{j=1}^{m-2}(-1)^{j+1}\Psi_1\cdots\Psi_j((\Delta\tens\Ad)_j-\Ad_{j+1})\Psi_{j+1}\cdots\Psi_{m-2}\\
&=&\sum_{j=1}^{m-2}(-1)^{j+1}\Psi_1\cdots\Psi_{j-1}((\Delta_{j+1}-\Ad_j)\Psi_j-\Psi_{j}\Ad_{j+1})\Psi_{j+1}\cdots\Psi_{m-2}\\
&=&(-1)^m\Psi_1\cdots\Psi_{m-2}\Ad_{m-1}-(\Ad_1+\sum_{j=2}^{m-1}(-1)^{j+1}\Delta_j )\Psi_1\cdots\Psi_{m-2}
\end{eqnarray*}
where the $\Ad$ terms cancel between the sum and the displaced sum except for the top term of one sum and the bottom term of the other. In the $\Delta$ sum 
all the indices of $\Psi$ are two or more smaller than the index of $\Delta$ so commute to the right. Combining with our previous calculation, we have
\begin{eqnarray*} [m,-\Psi]\del_{m-1}&=& (\Ad_1-\Delta_1+\del_{m-1})[m-2,-\Psi]\\
&&+(\Ad_1-\Delta_2+\Delta_3+\cdots(-1)^{m}\Delta_{m-1})(-1)^{m-2}\Psi_1\cdots\Psi_{m-2}\end{eqnarray*}
which proves  (\ref{PsiDeltarelns}). 

Next we use this result as initial base for  induction on $i$ in a formula
\begin{eqnarray}\label{PsiDeltai} 
&&\nquad [m-i+1]\cdots [m]\del_{m-1} =(\sum_{j=1}^{i-1}(-1)^{j+1}\Ad_1^{\tens j} \nonumber \\
&&+\left((-1)^{i+1}A_1^{\tens i}+\sum_{j=i+1}^{m-1}(-1)^{j+1}\Delta_j\right)[m-i])[m-i+1]\cdots[m-1] 
\end{eqnarray}
where $[m]\equiv[m,\Psi]$ for brevity and the nesting is rightmost as for braided factorials so $[m-1]\equiv \id\tens[m-1,-\Psi]$. The case $i=1$ is (\ref{PsiDeltarelns}) which we have already proven while the case $i=m-1$ or $i=m$, suitably interpreted in the sense of absent sums or products when out of range, proves the lemma. We use identities
\[  \Psi_i\Ad_j^{\tens  i-j}=\Ad_j^{\tens i}-\Ad_i,\quad \Psi_i\Ad_j^{\tens  k}=\Ad_j^{\tens k}\Psi_{i-1},\quad i>j+k\]
where the commutation relation is due to acting in different spaces, with renumbering due to the notation. The first equation is a direct computation. One also has
\[ \Psi_i\Ad_j^{\tens  k}=\Ad_j^{\tens k}\Psi_i,\quad i<j-1,\quad  \Psi_i\Ad_j^{\tens  k}=\Ad_j^{\tens k}\Psi_i,\quad j\le i<j+k-1\]
which we do not need right now, in the first case due to different tensor products and in the second case because $\Psi$ is a morphism in the crossed module category and hence commutes with $\Ad$ applied to tensor powers that include those on which $\Psi$ acts. Assuming (\ref{PsiDeltai}) for $i-1$ in the role of $i$, what we need to show to prove (\ref{PsiDeltai}) for $i$ is
\begin{eqnarray*}
&&[m-i+1]\left(\sum_{j=1}^{i-2}(-1)^{j+1}\Ad_1^{\tens j}+\left((-1)^{i}A_1^{\tens i-1}+\sum_{j=i}^{m-1}(-1)^{j+1}\Delta_j\right)[m-i+1]\right)\\
&&= \left(\sum_{j=1}^{i-1}(-1)^{j+1}\Ad_1^{\tens j}+\left((-1)^{i+1}A_1^{\tens i}+\sum_{j=i+1}^{m-1}(-1)^{j+1}\Delta_j\right)[m-i]\right)[m-i+1]. 
\end{eqnarray*}
Now, the first sum commutes with $[m-i+1]$ since on the left this is $1-\Psi_i+\Psi_i\Psi_{i+1}+\cdots+(-1)^{m-i}\Psi_i\cdots\Psi_{m-1}$ due to the right-most embedding. These commute past the $\Ad_1^{\tens j}$ getting changed to $[m-i+1]$ embedded on the right (where the numbering is reduced by one). Hence the first term on the left is  
$ \sum_{j=1}^{i-2}(-1)^{j+1}A_1^{\tens j}[m-i+1]$. Next $\Psi_i\Psi_{i+1}\cdots\Psi_{m-1}\Ad_1^{\tens i-1}=(\Ad_1^{\tens i}-\Ad_i)\Psi_i\Psi_{i+1}\cdots\Psi_{m-2}$ as the $\Psi_{i+1}$ and higher commute, reducing index by 1, while $\Psi_i$ computes as shown. Hence the middle terms gives  
\[[m-i+1](-1)^i\Ad_1^{\tens i-1}[m-i+1]=(-1)^i\Ad_1^{\tens i-1}[m-i+1]+(-1)^{i+1}(\Ad_1^{\tens i}-\Ad_i)[m-i][m-i+1],\]
the first term of which completes our previous sum to give the first desired term. Accordingly we need only show for the remaining term that
\[ [m-i+1]\sum_{j=i}^{m-1}(-1)^{j+1}\Delta_j=\left((-1)^{i+1}\Ad_i+ \sum_{j=i+1}^{m-1}(-1)^{j+1}\Delta_j\right)[m-i].\]
But this is just the same identity (\ref{PsiDeltarelns}) already proven but for  $[r]\del_{r-1}$, i.e.\ $r=m-i+1$ in the role of $m$, for the $m$ tensor factors numbered $i,\cdots,m$. This completes our proof of (\ref{PsiDeltai}) for all $i$ and proves the lemma. \endproof

This fleshes out the braided-Hopf algebra interpretation of the Woronowicz exterior algebra on a Hopf algebra\cite{Wor} using  \cite{MaTao} for the direct treatment of $\extd$.

\subsection{Differential calculi on braided linear spaces} 

For completeness, we give another braided construction namely the exterior algebra $\Omega(B)$ on a Hopf algebras $B$ in a braided Abelian category. This includes braided symmetric algebras $B=B_+(\Lambda^1)$ as above generated canonically by an object $\Lambda^1$. The further data we will need is a surjective morphism $\varpi:B\to \Lambda^1$ in the category such that
\begin{equation}\label{ederiv} \varpi\circ\cdot=\eps\tens\varpi+\varpi\tens\eps\end{equation}
This data arises naturally as follows: suppose $B^\flat$ is a (possibly degenerately) dually paired braided-Hopf algebra from the right (so the pairing is $\ev:B\tens B^\flat\to \und 1$) and  $\CL$ a rigid primitive sub-object $\CL\subset B^\flat$ (so that the coproduct restricted to $\CL$ is the additive one). We view the duality pairing restricted to a map $B\tens\CL\to \und 1$ as a map
\[ \varpi:B\to \Lambda^1=\CL^*,\quad \varpi=(\ev\tens\id)(\id\tens\coev_{\CL})\]
which then obeys (\ref{ederiv}). This is surjective if there does not exist $\eta\in \CL$ which pairs to zero with all of $B$. In the case of $B=B_+(\Lambda^1)$ or any other graded braided Hopf algebra of the form $B=\und 1\oplus \Lambda^1\oplus B_{>1}$ generated in degree 1 by an object $\Lambda^1$, we simply take $\varpi:B\to \Lambda^1$ as the projection to degree 1.

\begin{proposition}\label{diffplaprop} Let $B$ be a Hopf algebra in an Abelian braided category and $\varpi:B\to \Lambda^1$ a surjective morphism obeying (\ref{ederiv}). Then
\[ \Omega= B\und\tens\Lambda,\quad  \Lambda= T\Lambda^1/\<{\rm image}(\id+\Psi_{\Lambda^1,\Lambda^1})\>\]
\[ \extd|_B=(\id\tens\varpi)\Delta,\quad \extd|_\Lambda=0\]
is a differential exterior algebra on $B$ in the category (one in which all structure maps are morphisms). 
\end{proposition}
\begin{figure}\label{diffplane}\[ \includegraphics[scale=0.5]{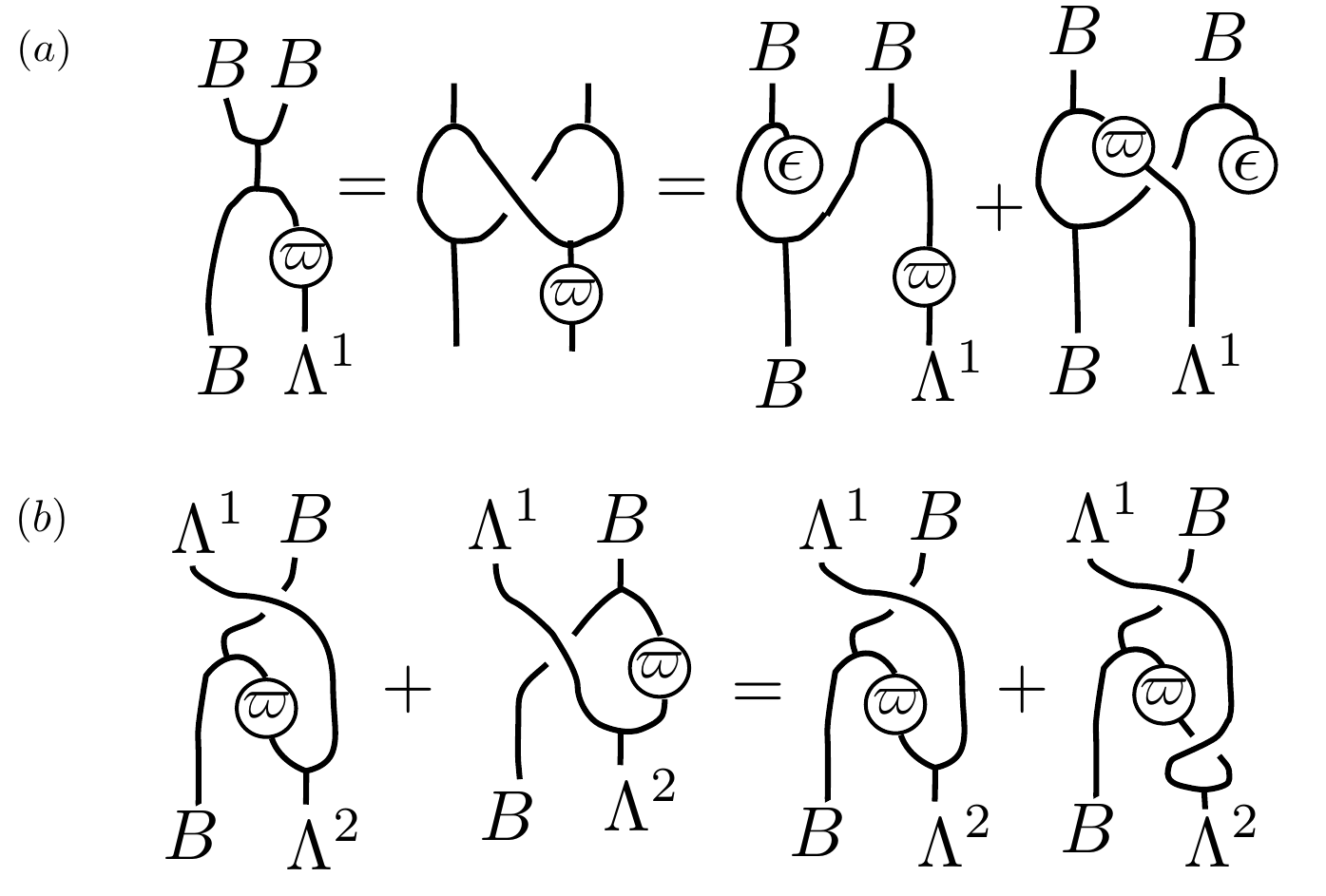}\] \caption{Diagrams in the proof of Proposition~\ref{diffplaprop} for quantum differentials on braided planes}
\end{figure}
\proof  The proof is done diagrammatically in Figure~\ref{diffplane} and applies generally but for convenience of exposition we also refer to concrete elements. The braided tensor product $B\und\tens\Lambda^1$ in concrete terms means $(b\tens v)(c\tens w)=b\Psi(v\tens c)w$ and featured already in the definitiion of a braided-Hopf algebra. Part  (a) computes $\extd(bc)$ using the braided coproduct homomorphism property and (\ref{ederiv}). Using the counit axioms and $\varpi$ a morphism we obtain  $b\extd c$ for the first term and $(\extd b)c $ for the second when we remember the braided tensor product.  Part (b) checks that $\extd$ extends as a graded derivation with respect the braided tensor product. We compute 
\[ \extd(\omega b)+\omega\extd b=( \id\tens\cdot)(\id+\Psi^{-1})({\rm something})\]
hence vanishes in $B\tens\Lambda^2$ which agrees with $\extd\omega=0$ for all $\omega\in \Lambda^1$. The same applies starting in $\Lambda\tens B$ \endproof

This gives a differential structure on our braided-symmetric algebras $B_+(V)$ regarded as noncommutative spaces. If the category is the comodules of a coquastriangular Hopf algebra, for examples, our construction is covariant in that all structure maps are comodule maps. Also note that if $\{e_i\}$ is a basis of $V$, we have explicitly
\[ \extd = \del^i (\ )e_i\]
where $\del^i$ are the  (right handed) {\em braided partial derivatives} defined by 
\[ \Delta b= 1\tens b+ \del^i b\tens e_i+\cdots \]
They are given explicitly at the level of the tensor algebra by 
\[ \extd(v_1\tens ...\tens v_n)=(\eta_{n-1}\tens\eta_1)\left[{n\atop n-1},\Psi\right](v_1\tens\cdots\tens v_n)\]
\[ \left[{n\atop n-1},\Psi\right]=1+\Psi_{n-1}+\Psi_{n-1}\Psi_{n-2}+\cdots+\Psi_{n-1}\cdots\Psi_1\]
where the last tensor factor of the result is viewed in $\Lambda^1$.  

\begin{remark}\label{remplane} If we take the quadratic version $S(V):=B_+^{quad}(V)=TV/\<\ker(\id+\Psi)\>$ then our above construction gives
\[ \Omega(S(V))=S(V)\und\tens S(V^*)^!,\quad \extd v=1\tens v\]
where $!$ denotes the Koszul dual. If a quadratic algebra on a vector space $W$  has relations $R\subset W\tens W$ as the subspace being set to to zero then its Koszul dual is the quadratic algebra on $W^\flat$ with relations $R^\perp\subset W^\flat\tens W^\flat$. This is normally  done in the category of vector spaces but we do it here in a braided category using the right dual so that $W^\flat=V$. \end{remark}

\begin{example}\label{fermqplane} We let $B=\A_q^2=B_+(V)$ be the quantum plane associated to the standard corepresentation $V={\rm span}\{x,y\}$ in the braided category of right $k_q[GL_2]$-comodules with $q^2\ne 1$. Here 
\[ \Psi(x\tens x)= q^2 x\tens x,\  \Psi(x\tens y)=q y\tens x,\  \Psi(y\tens x)=q x\tens y+(q^2-1)y\tens x,\  \Psi(y\tens y)=q^2y\tens y\]
is given by a particular non-standard normalisation of the usual $\CR$ on $k_q[GL_2]$ (one that does not descend to $k_q[SL_2]$). The kernel of $\id+\Psi$ gives us the relations $yx=qxy$ of the quantum plane since  $(\id+\Psi)(y\tens x-q x\tens y)=0$. The algebra $\Lambda=\A_q^{0|2}=\Lambda(V)$ is the fermionic quantum plane 
\[ \extd x\wedge \extd x=0,\quad \extd y\wedge\extd y=0,\quad \extd y\wedge\extd x=-q^{-1}\extd x\wedge\extd y\]
 where this time the same basis is denoted $\{\extd x,\extd y\}$ as a basis of $\Lambda^1=V$ and one can check for example that $(\id+\Psi)(\extd x\tens\extd y)=\extd x\tens \extd x+q\extd y\tens\extd x$ from the stated braiding. Indeed, it known that $\A_q^2$ and $\A_q^{0|2}$ as Koszul dual as first pointed out by Manin\cite{Man}. The differential on $v\in V\subset B$ is $\extd v=1\tens v$ i.e., $v$ viewed in the $\Lambda^1$ copy of $V\subset\Lambda(V)$. The relations between $B$ and $\Lambda(V)$ are the braided tensor product so $(1\tens v)(w\tens 1)=\Psi(v\tens w)$. This again comes from the same braiding as above but viewed now as defining the relations
\[ (\extd x)x=q^2x\extd x,\quad (\extd x)y=q y\extd x,\quad  (\extd y)x=qx\extd y+(q^2-1)y\extd x,\quad  (\extd y)y=q^2y\extd y.\]
One can check for example that $\extd(yx-qxy)=0$ as it should.  By construction, this exterior algebra on the quantum plane is $k_q[GL_2]$-covariant. The associated partial derivatives are 
\[ \del^1(x^my^n)=[m,q^2]x^{m-1}y^n q^n,\quad \del^2(x^my^n)=x^m [n,q^2]y^{n-1}\]
using the braided coproduct on general monomials computed in \cite{Ma:book} from the braided-integers $[n,\Psi]$. The partial derivatives here were first found by Wess and Zumino in another approach. They are naturally `braided right derivations' with an extra $q^n$ in the first expression, in order that $\extd$ is a left derivation, acting as braided $q^2$-derivatives in each variable. One can check that $\del^2\del^1=q\del^1\del^2$ as operators, also as per the general theory in \cite{Ma:book}. 
\end{example}

This reworks the treatment of quantum-braided planes and their differentials in \cite[Sec. 10.4]{Ma:book}\cite{Ma:fre} now as an example of our above  canonical construction based on $B_+(\Lambda^1)$, as opposed to a compatible pair of R-matrices $R,R'$ as previously.

\section{Braided fourier transform and application to Hodge theory}

Fourier transform on Hopf algebras is part of their classical literature. It  was extended to braided-Hopf algebras in \cite{LyuMa} and related works and applied to braided linear spaces in \cite{KemMa}, though not the ones we consider here. We used diagram proofs and will do so again, while another work from that era is \cite{Lyu}. We first explain the general (super) formulation and then apply it to the Hodge operator, including $k(S_3)$ as an example. 

\subsection{Super-braided Fourier theory}

 In any braided category $\CC$ and $B\in \CC$ a braided Hopf algebra dually paired with a braided Hopf algebra $B^\star$, we have three actions which we will consider and which we collect in Figure~1 in a diagrammatic notation\cite{Ma:cro,Ma:bhop}. Diagrams are read as operations flowing down the page, with tensor products and the unit object $\und 1$ suppressed. Two strands flowing own and merging denotes the product and one strand flowing down and splitting denotes the coproduct. As in Section~2 when discussing duals, we assume a pairing $\ev:B^\star\tens B\to \und 1$ which we can write diagrammatically as $\cup$ and with respect to which the product on one side is adjoint to the coproduct on the other. The counit of $B$ is also adjoint to the unit $\eta:\und 1\to B^*$ in the sense $\ev(\eta\tens ( ))=\eps_B$ and vice versa $\ev(( )\tens \eta)=\eps_{B^\star}$, where we have canonical isomorphisms $\und 1\tens B\isom B$ etc which we use. In  a concrete $k$-linear setting we can suppose that $\und 1=k$ and $\eta(1)=1$ to simplify the above. {\rm Reg} makes $B$ a right $B^\star$ module algebra in the braided category. The principal ingredient of Reg here is actually a left action $\vdash$ 
making $B$ a left $B^\star$-module algebra in the braided category. Similarly, we have a straightforward right action $\dashv$ under which $B^\star$ is a right $B$-module algebra \cite{Ma:bhop}. 

\begin{figure}
\[ \includegraphics[scale=.6]{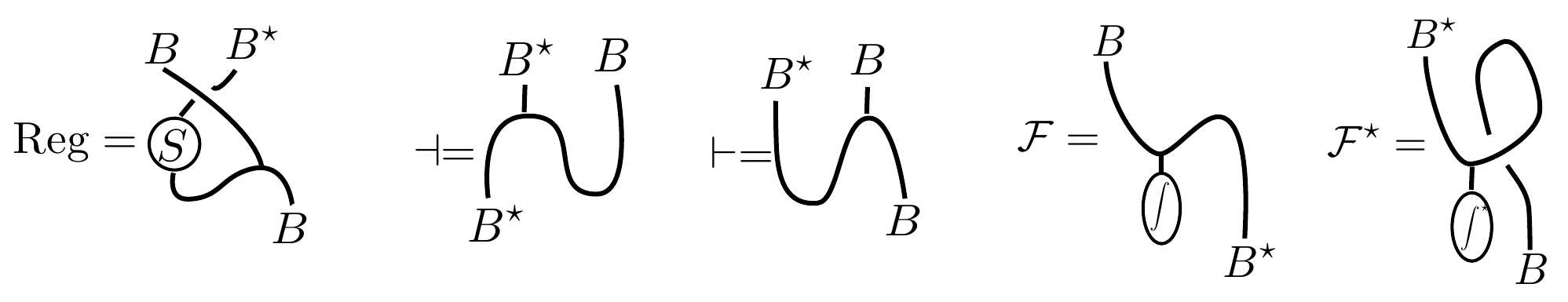}\]
\caption{Diagrammatic definitions of relevant actions, Fourier transform $\CF$ and adjoint Fourier transform $\CF^*$ on a braided-Hopf algebra $B$ with right dual $B^\star$.}
\end{figure}

We also need the notion of a {\em left integral}  and the simplest thing is to require a morphism $\int: B\to \und 1$ in the sense $(\id\tens\int)\Delta=\eta\tens\int$. However, we do not want to be too strict about this. For example, for the finite anyonic braided line $B=k[x]/(x^{n+1})$ in the braided category of $\Z/(n+1)$-graded spaces with braiding given by an $n+1$-th root of 1 and $|x|=1$, the obvious  $\int x^m=\delta_{m,n}$ is not a morphism to $\und 1$. Our approach is to live with this and not necessarily assume any morphism properties; we can still use the diagrammatic notation but be careful not to pull the map through any braid crossings. A more formal approach is to view it as a morphism $B\to K$ where $K=k$ taken with degree $n$ in the case of the anyonic braided line. The uniqueness of the integral when it exists is similar to the Hopf algebra case (see \cite{BKLT} for a formal proof).

For Fourier transform we need not only that $B^\star$ is dually paired but that $B$ is actually rigid with dual object $B^*$. Again, this is a very strong assumption, analogous to finite-dimensionality of $B$ and amounting to this in the typical $k$-linear case. It means that there is a coevaluation map $\exp={\rm coev}: \und 1\to B\tens B^*$, denoted by $\cap$ in the diagrammatic notion, which obeys the well-known `bend-straightenning axioms' with respect to $\cup$. We similarly require a right integral which is not necessarily a morphism $\int^\star: B^\star\to \und 1$. We can live with this or suppose formally that $\int^\star:B^\star\to K^*$ where $K^*\tens K=\und 1= K\tens K^*$ as objects. In our $k$-linear setting this will be by the identification with $k$. The theory below could be generalised to include some infinite-dimensional cases or else these could be treated formally eg in a graded case with $B^*$ a graded dual, each component rigid and the result a formal power series in a grading parameter.

\begin{definition}\cite{LyuMa,KemMa} Let $B$ be a Hopf algebra in a braided category with $B$ rigid and $\int$ a left integral as above. Then $\CF:B\to B^\star$ defined in Figure~1 is called the {\em braided Fourier transform}. We similarly define a dual Fourier transform $\CF^\star:B^\star\to B$ if $B^\star$ has a right integral. 
\end{definition}
These maps are no longer morphisms if the integrals are not, or one can say more formally that $\CF:B\to K\tens B^\star$ and $\CF^\star: B^\star\to K^*\tens B^\star$. The following extends and completes \cite{KemMa}.

\begin{proposition} \label{Foubra} In the setting of the definitions above 
\[ \CF\circ {\rm Reg}= \cdot \circ(\CF\tens\id),  \quad  \dashv (\CF\tens\id)=\CF\circ\cdot\]
Moreover, if $\int^*$ is a {\em right} integral on $B^\star$ then
\[ \CF^\star\CF=\mu S,\quad \mu:=(\int\tens\int^*)\exp\]
If the integrals are both unimodular and morphisms then $\CF\CF^\star=\mu S$ and $[\CF,S]=0$ when $\mu$ is invertible (see Figure~2 for the general case).
\end{proposition}
\proof Here $\CF {\rm Reg}$ and $\CF^\star\CF$ are already covered in \cite[Figs~3-5]{KemMa} so we do not repeat all the details here. We recall only the diagram proof for $\CF^*\CF$ using the lemma in \cite[Fig.~2(b)]{KemMa} at the first equality in Figure~2 and note that we did not need to assume that $\int,\int^*$ are morphisms to $\und 1$ as in \cite{KemMa} as long as we keep the integrals to the left. The second line now uses the same lemma but this time on $B^\star$ to compute $\CF\CF^*$ as shown provided $\int^*$ is also a left integral so that the lemma applies and $\int$ is also a right integral. If $\int,\int^*$ are morphisms to $\und 1$ so we can take them through braid crossings to obtain $\mu S$  and then $\mu \CF S=\CF\CF^\star\CF=\mu S\CF$. The general result $\dashv (\CF\tens\id)= \CF\circ\cdot$ follows more simply from the duality pairing and associativity of the product of $B$. \endproof
\begin{figure} 
\[ \includegraphics[scale=.45]{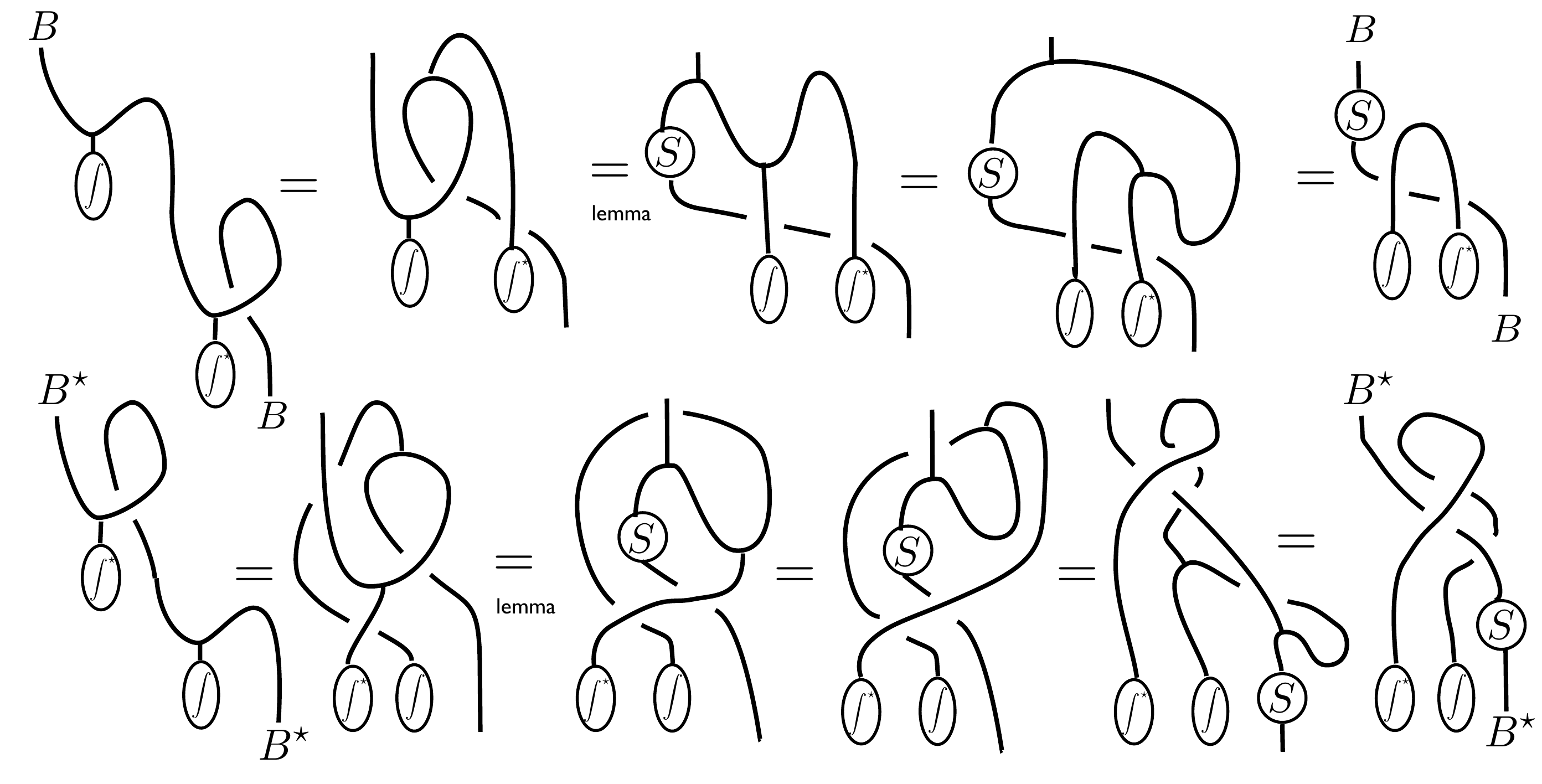}\]
\caption{Diagrammatic computations of  $\CF^\star\CF$ and $\CF\CF^\star$ in Proposition~\ref{Foubra}}
\end{figure}

The map $\CF^\star$ here is a right-integral version of the theory which is being used to define the adjoint Fourier transform and converted to a left version via $\Psi$. The braided antipode $S$ plays the role of the minus sign familiar in classical Fourier theory and $\mu$ plays the role of $2\pi$. If $\mu$ and $S$ are invertible then the stated results imply that $\CF$ is invertible at least in the $k$-linear setting (with $\CF^{-1}=S^{-1}\CF^\star$ in the unimodular trivial morphism case). Also, if we compose $\CF$ with $S$ then the first property above becomes 
\begin{equation}\label{SFvdash} S\CF\circ\vdash =\cdot(\id\tens S\CF).\end{equation}

\begin{example} For the elementary example of $B=k[x]/(x^{n+1})$ we have ${\rm Vol}=x^n$ as discussed, $B^\star=k[y]/(y^{n+1})$ where $|y|=-1$ and  ${\rm Vol}^*=y^n$ for the top form and $\exp=\sum_{m=0}^n x^m\tens y^m /[m,q]!$ as mentioned at the end of Section~2.2. Here $[m,q]=(1-q^m)/(1-q)$ and $q$ is a primitive $n+1$-th root of 1. We have
\[ \CF(x^m)={y^{n-m}\over [n-m,q]!},\quad \CF^\star(y^m)={q^{(n-m)^2} x^{n-m}\over [n-m,q]!},\quad \mu=[n,q]!^{-1}\]
as well as  $S x^m=(-x)^m q^{m(m-1)\over 2}$ and ditto with $x$ replaced by $y$. One can see that $\CF^\star\CF=\mu S$ using $q^{n(n+1)\over 2}=(-1)^n$.  Due to the nontrivial braidings of the integrals, however, the right hand side in Figure~2 gives 
\[ \CF\CF^\star=q^{2 D+1}\mu S\]
where $D$ is the monomial degree operator. The same method as in the proof above now gives us $\CF S=q^{2D+1}S\CF$ or equivalently $S\CF=\CF S q^{2D+1}$, which one may verify from the stated $\CF,S$. 
\end{example}

\begin{example}\label{Fouplane} The fermionic quantum plane in Example~\ref{fermqplane} is a super-braided Hopf algebra, i.e. we take $B=\A_q^{0|2}=B_-(\Lambda^1)$ in the category of $k_q[GL_2]$-comodules. If we now denote the braiding in Example~\ref{fermqplane} as $\Psi_+$ adapted to $B_+$ as a braided-Hopf algebra,  we now take a different normalization $\Psi=\Psi_-=q^{-2}\Psi_+$ (i.e. induced by another differently normalised coquasitriangular structure) so that $\ker(\id-\Psi_-)={\rm image}(\id+\Psi_+)$. The $\Psi_\pm$ here are $q^{\pm 1}$ times the braiding in the standard $q$-Hecke normalisation. There are no new relations in higher degree so $B_-(\Lambda^1)=B_-^{quad}(\Lambda^1)$. For brevity we let $e_1=\extd x$, $e_2=\extd y$ then  $\Psi(e_2\tens e_1)=q^{-1}e_1\tens e_2+\lambda e_2\tens e_1$ etc. in the new normalisation, 
where $\lambda=1-q^{-2}$. We now develop $\A_q^{0|2}$  as a super-braided Hopf algebra with $e_i$ primitive and underlying braiding $\Psi$ (meaning we actually transpose with super-braiding $\Psi_{sup}$ having additional $\pm$ factors according to the monomial degrees). This implies $S(e_1e_2)=q^{-2}e_1e_2$ as well as $S(1)=1$, $S(e_i)=-e_i$. On the dual side we have $B^\star=B_-(\Lambda^1{}^*)$ with a dual basis of generators $f^1,f^2$, underlying braiding
\[ \Psi(f^i\tens f^i)=f^i\tens f^i,\quad \Psi(f^1\tens f^2)=q^{-1}f^2\tens f^1+\lambda f^1\tens f^2,\quad \Psi(f^2\tens f^1)=q^{-1}f^1\tens f^2,\]
relations $f^2 f^1=-q f^1 f^2$ and $S(f^1f^2)=q^{-2}f_1f_2$. There is up to scale a unique top degree in each case, namely ${\rm Vol}=e_1e_2$ and ${\rm Vol}^*=f^1f^2$  and we find $\<{\rm Vol}^*,{\rm Vol}\>=\ev(f^1 \tens f^2,[2,-\Psi](e_1\tens e_2)\>=-q^{-1}$, so that 
\[ \exp=1\tens 1+\sum_{i=1}^2 e_i\tens f^i - q{\rm Vol}\tens{\rm Vol}^*\]
We define integrals via $\int {\rm Vol}=1$ and $\int^*{\rm Vol}^*=1$ but note that these are not morphisms. Rather we use braidings
\[ \Psi(f^1\tens e_1)=e_1\tens f^1+(1-q^2)e_2\tens f^2,\quad \Psi(f^2\tens e_2)=e_2\tens f^2,
\quad \Psi(f^i\tens e_j)=qe_j\tens f^i\]
for $i\ne j$ (these are obtained from the 2nd inverse $\tilde R$ as in \cite[Propn.~10.3.6]{Ma:book} for $R$ normalised to our case) to find $\Psi(f^i\tens {\rm Vol})=q {\rm Vol}\tens f^i$ and hence $\Psi({\rm Vol}^*\tens{\rm Vol})=q^2 {\rm Vol}\tens{\rm Vol}^*$. We similarly have $\Psi(f^i\tens {\rm Vol}^*)=q{\rm Vol}^*\tens f^i$. From these it is clear that $\int$ and $\int^*$ are not morphisms in the underlying comodule category. Again, there can be further signs according to the super degrees for the actual super-braiding $\Psi_{sup}$ when we read diagrams in the super-braided case. In particular, we find
\[ \Psi_{sup}^{-1}\exp=1\tens 1- f^1\tens e_1-q^2 f^2\tens e_2-q^{-1}{\rm Vol}^*\tens{\rm Vol}\]
needed in the computation of $\CF^\star$. We now read off from the diagrammatic definitions in Figure~1,
\[ \CF\begin{cases}1\cr e_1\cr e_2\cr {\rm Vol}\end{cases}=\begin{cases}-q{\rm Vol}^*\cr f^2\cr -q^{-1}f^1\cr 1\end{cases},\quad 
\CF^\star\begin{cases}1\cr f^1\cr f^2\cr {\rm Vol}^*\end{cases}=\begin{cases}-q^{-1}{\rm Vol}\cr -q^2 e_2 \cr q e_1\cr 1\end{cases},\quad\mu=-q.\] 
One can verify that $\CF^\star\CF=\mu S$ as it must by Proposition~\ref{Foubra}. We also have
\[ \CF\CF^\star=\mu q^{2(D-1)}S\]
where $D$ is the monomial degree and one can check that this agrees with  the lower line in Figure~2 where we use the above computations to read off the right hand side. In this case $\CF S=q^{2(D-1)}S\CF$ or $S\CF=\CF S q^{2(D-1)}$ as one can verify from the stated form of $\CF,S$. \end{example}

A similar approach can be used for other quantum planes to express their differential exterior algebras as super-braided Hopf algebras with possibly a different underlying coquastriangular structure from one used for the coordinate algebra as a braided-Hopf algebra. In the presence of an invariant quantum metric we reproduce the otherwise ad-hoc approach to $q$-epsilon tensors and Hodge theory on braided-quantum planes in \cite{Ma:book}. The $f^i$ generate antisymmetric vectors and $\dashv, \vdash$ define an interior product connected to the exterior algebra product via $\CF$. This example should be seen as a warm-up to Section~3.2 where we look at bicovariant differentials on Hopf algebras themselves. As illustrated here, the actual theory is read off the diagrams with the appropriate braiding including signs. We could indeed shift all constructions to this new super-braided category and say that the above example is an ordinary braided-Hopf algebra there, but we not do so since there will normally be other (bosonic) objects also of interest in the original category. In our context the nicest case is where $\int,\int^*$ are morphisms to $\und 1$ when viewed in the original category but do not necessarily respect the super-degree, for example they could be odd maps in the super-sense in which case they are not morphisms in the super version with extra signs (so we need the slightly more general picture as above). We assume they have the same parity of support (both odd or both even maps). Then $\CF,\CF^\star$ also have this party. 

\begin{corollary}\label{Fousup} If $\int,\int^*$ are unimodular, morphisms in the underlying category and of the same parity $p$ then $\CF\CF^\star=(-1)^p\mu S$. When $\mu,S$ are invertible we have $\CF S=(-1)^p S\CF$ and $\CF^\star=\mu S \CF^{-1}$
\end{corollary}
\proof  For $\CF\CF^\star$ we have to compute the right hand side of the lower diagram in Figure~2, which now has extra signs. We can still bring out $\mu'=(\int^*\tens\int)\Psi_{sup}^{-1}\exp$ since any signs from crossing the $\int$ leg cancel with signs from crossing the $\int^\star$ leg by our assumptions. Next, we can lift $\int^\star$ through the crossing at the price of $(-1)^p$ in computing $\mu'$.  We already have $\CF^\star\CF=\mu S$ from Proposition~\ref{Foubra} and can then conclude the rest. \endproof

The behaviour of $\CF$ with respect to ${\rm Reg},\dashv$ has an unchanged form as these statements do not involve additional transpositions, except that the actions  themselves are computed for the super-braided Hopf algebra eg with the super-braided coproduct and hence the super-braided Leibniz rule expressed in super-braided module algebra structures. The property (\ref{SFvdash})  becomes 
\begin{equation}\label{superSFvdash} S\CF\circ\vdash=\cdot((-1)^{pD}\tens S\CF)\end{equation}
due to the crossing of the first input on the right hand side with the integral in $\CF$.

\subsection{Hodge theory on Hopf algebras}

We are now going to  compute our super-Fourier theory for $B=\Lambda_{min}=B_-(\Lambda^1)$ where $\Lambda^1$ is a rigid object in the braided category of right $A$-crossed modules and $A$ is a Hopf algebra with invertible antipode. Here $B$ is a super-braided Hopf algebra in the category and we assume it has a top degree component $K$ of dimension 1, i.e. up to scale a unique top form ${\rm Vol}\in B_-(\Lambda^1)$. This gives us a unimodular integral $B\to k$ by  $\int {\rm Vol}=1$ and zero for lower degrees. To see this, note that the formula in Proposition~\ref{TVcoprod} ensures that $\und\Delta {\rm Vol}={\rm Vol}\tens 1+1\tens{\rm Vol}$ plus terms of intermediate degree, and we never reach the top degree on applying $\und\Delta$ to lower degree. One can think of this more formally as a morphism $B\to K$ with some possibly non-trivial generator. We also have an identification $B^\star=B_-(\Lambda^1{}^*)$ by extending the duality pairing $\Lambda^1{}^*\tens\Lambda^1$ as a braided-Hopf algebra pairing, given that this is now non-degenerate after quotienting by the relations of $B_-$ as explained in Section~2.2. 
Hence we obtain a unimodular integral on this too. In the nicest case, the top forms ${\rm Vol},{\rm Vol}^*$ of degree $n$ (say) span the trivial object $\und 1$ so that $\int,\int^*$ are morphisms to $\und 1$ but of parity $n$ mod 2, so we are  in the setting of Corollary~\ref{Fousup}. 

Next, in non-commutative geometry a metric is $g\in \Omega^1\tens_A\Omega^1$ with an inverse $(\ ,\ ):\Omega^1\tens_A\Omega^1\to A$. One can show that in this case $g$ must be central. Normally, one also requires the metric to be `quantum symmetric' in the sense of the product $\wedge(g)=0$ in $\Omega^2$. We are interested in left-invariant metrics where $g\in \Lambda^1\tens\Lambda^1$. 

\begin{lemma} A bi-invariant metric on a Hopf algebra $A$ with bicovariant calculus is equivalent to an $A$-crossed module isomorphism $g:\Lambda^{1*}\isom \Lambda^1$.  The metric is quantum symmetric if and only if $\Psi(g)=g$. 
\end{lemma}
\proof The metric being bi-invariant means that it is an element $g\in \Lambda^1\tens\Lambda^1$ which is invariant  under the coaction $\Delta_R$ on the tensor product. The existence of a bimodule map $(\ ,\ )$ requires $g$ to be central which in turn requires that $g$ is invariant under the crossed module right action $\ra$ (since this determines the cross product of $A\rbiprod \Lambda$). So a metric is equivalent to a morphism $\und 1\to \Lambda^1\tens\Lambda^1$ in the crossed module category. Evaluation from the left makes this equivalent a morphism as stated, which we also denote $g$. Here $\Lambda^*$ is again a right crossed module in the usual way (via the antipode). Clearly $\wedge(g)=0$ if and only if $g\in \ker[2,-\Psi]=\ker(\id-\Psi)$ according to the relations of $B_-(\Lambda^1)$. \endproof

Given a bi-invariant metric we therefore have $B_\pm(\Lambda^1{}^*)\isom B_\pm(\Lambda^1)$  hence combined with the above remarks in the finite-dimensional case, an isomorphism  which we also denote $g:B_\pm(\Lambda^1)^*\to B_\pm(\Lambda^1)$. We are now ready to define the Hodge operator, using the $B_-$ version. We do it in the nicest case but the same ideas can be used more generally as we have seen in Section~3.1.

\begin{definition}\label{sharp} Suppose that $\Lambda^1$ is finite-dimensional in the category of right $A$-crossed modules, $g$ a bi-invariant metric and $B_-(\Lambda^1)$ finite-dimensional with a 1-dimensional top degree $n$ and central bi-invariant top form ${\rm Vol}$ used to define $\int$.  We define the Hodge star 
\[ \sharp{\ }=g\circ\CF: B_-(\Lambda^1)^m\to B_-(\Lambda^1)^{n-m}\]
which we extend as a bimdodule map to $\Omega^m\to \Omega^{n-m}$. 
\end{definition}
By construction our $\sharp$ is a morphism in the crossed-module category. In geometric terms this means that it extends as a bimodule map and is bicovariant under the quantum group action on $\Omega$. We also define $\sharp^\star=(-1)^D\circ \CF^\star\circ g^{-1}$ where $D$ is the degree operator. 

\begin{proposition}\label{outer} In the setting of Definition~\ref{sharp},  $\mu=\<{\rm Vol},{\rm Vol}\>^{-1}\in k^\times$, $\sharp$ is invertible and $\sharp S=(-1)^nS\sharp$. If the metric $g$ is quantum symmetric  then $S=(-1)^D$, $\sharp^\star=\sharp$ and $\sharp^2=\mu$ on degrees $D=0,1,n-1,n$. \end{proposition}
\proof  Here $\<{\rm Vol},{\rm Vol}\>$ is non-zero since otherwise ${\rm Vol}$ would be zero in $B_-(\Lambda^1)$, and its inverse supplies the coefficient of the top component of $\exp$, which is $\mu$.  Since $\mu\ne 0$ we can apply Corollary~\ref{Fousup} to see in particular that $\sharp,S$ graded-commute. That  $S|_{0,1,n-1,n}=(-1)^D$ 
i.e. on the outer degrees is clear on degrees 0,1 and then holds on degrees $n,n-1$ due to $\sharp, S$ graded-commuting. Next, in terms of $\sharp^\star$ with the metric identification, the result in Corollary~\ref{Fousup} becomes $\sharp^\star\sharp=\mu (-1)^D S$ and $\sharp(-1)^D\sharp^\star=\mu (-1)^n S$ since the parity of the integral is $n$ mod 2. Taking the $(-1)^D$ to the left in the latter equation makes it $(-1)^{n-D}$ so that 
\[ \sharp^\star=\mu (-1)^D S\sharp^{-1}\]
on all degrees, giving $\sharp^\star|_{0,1,n-1,n}=\mu\sharp^{-1}$ on the outer degrees. On the other hand, we have
 \[\exp=1\tens 1 + g+ \cdots + g^{(n-1)}+ \mu {\rm Vol}\tens{\rm Vol}\]
(for some element $g^{(n-1)}\in \Lambda^{n-1}\tens\Lambda^{n-1}$), while the definition of $\sharp^\star$ is such that it is given by integration agains $\Psi^{-1}\exp$ without any signs. Since $g$ (by the quantum symmetry assumption) and $1\tens 1$ are invariant under $\Psi$, these terms are the same, and hence $\sharp^\star=\sharp$ on degrees $n-1,n$ and hence $\sharp^2=\mu$ on these degrees. In that case $\sharp^\star(\sharp\omega)=\mu\omega=\sharp^2\omega$ on all $\omega$ of degree $n-1,n$ tells us that $\sharp^\star=\sharp$ on degrees 0,1 also, and hence that $\sharp^2=\mu$ on these degrees also. This means that 
\[ \Psi_{sup}^{-1}\exp=1\tens 1-g +\cdots +(-1)^{n-1}g^{(n-1)}+(-1)^n\mu{\rm Vol}\tens{\rm Vol}\]
for the computation of $\CF^\star$ and similarly without the signs for $\sharp^\star$. 
\endproof

We similarly define left and right interior products 
\[ \vdash: \Lambda^1\tens B_-(\Lambda^1)^m\to B_-(\Lambda^1)^{m-1},\quad \dashv: B_-(\Lambda^1)^m\tens\Lambda^1\to B_-(\Lambda^1)^{m-1}\]
by restricting the left and right actions in Section~3.1 (these are the left and right braided-partial derivatives in the sense of \cite{Ma:fre,Ma:book}). We then extend these to bicovariant bimodule maps
\[\vdash:\Omega^1\tens_A\Omega^m\to \Omega^{m-1},\quad \dashv:\Omega^m\tens_A\Omega^1\to \Omega^{m-1}\]
given by
\[ (a\eta)\vdash (b\omega)=(a\eta,b\omega\Bo)\omega\Bt,\quad (b\omega)\dashv (a\eta)=b\omega\Bo(\omega\Bt,a\eta),\quad \forall a,b\in A,\  \eta\in \Lambda^1,\ \omega\in \Lambda\]
where we underline the braided-coproduct of $\Lambda$. In other words, we extend the braided coproduct as a bimodule map $\Omega\to \Omega\tens_A\Omega$ (not to be confused with the super-coproduct of $\Omega$ as a super-Hopf algebra) and then use the quantum metric  pairing to evaluate, taken as zero when degrees do not match. 

We can now interpret our Fourier theory in Section~3.1  as
\begin{equation}\label{inthodge}   S\sharp(\eta\vdash\omega)=\eta(S\sharp\omega),\quad \sharp(\omega\eta)=(\sharp\omega)\dashv\eta,\quad\forall\eta\in \Omega^1,\quad \omega\in \Omega\end{equation}
where $S$ is the super-braided antipode of $B_-(\Lambda^1)$ extended as a bimodule map to $\Omega$. It should not be confused with the super-coproduct of $\Omega$. We also define a left Lie derivative by
\[ \CL_\eta(\omega):=\eta\vdash \extd\omega+ \extd (\eta\vdash\omega),\quad\forall \eta\in \Omega^1,\ \omega\in\Omega\]
and associated codifferential and Hodge Laplacian
\begin{equation}\label{deltaHodgeL} \delta:=(S\sharp)^{-1}\extd(S\sharp),\quad \square:=\extd\delta+\delta\extd.\end{equation}
The use of $(S\sharp)^{-1}$ here is adapted to the left handed $\vdash$ and left-handed partial derivatives defined by $\extd f=\sum_a (\del^a f)e_a$ for any choice of basis $\{e_a\}$ of $\Lambda^1$. One could equally well use $\sharp$ but this would be adapted to $\dashv$ and right-handed partial derivatives. We also define the Leibnizator
\[ L_\delta(\omega,\eta)=\delta(\omega\eta)-(\delta\omega)\eta-(-1)^{|\omega|}\omega\delta\eta\]
as in \cite{Ma:rq}. 

\begin{corollary} \label{lbformula} The codifferential and Hodge Laplacian in (\ref{deltaHodgeL}) obey
\[ \delta(f\omega)=f\delta\omega+(\extd f)\vdash \omega,\quad \square(f\omega)=(\square f)\omega+f\square\omega+L_\delta(\extd f,\omega)+\CL_{\extd f}\omega\]
for all $f\in A$ and $\omega\in \Omega$. Moreover, 
\[  \delta \alpha=\alpha^a\delta e_a+g_{ab}\del^a\alpha^b,\quad \square f=(\del^af)\delta e_a+ g_{ab}\del^a\del^b f\]
\[ \square \alpha=\alpha^a \square e_a +(\square \alpha^a)e_a+ \del^a\alpha^b(L_\delta(e_a,e_b)+\CL_{e_a}e_b)+ \del^a\del^b\alpha^c \left((e_a\vdash(e_be_c))+ g_{bc}e_a-g_{ab}e_c\right)\]
where $\alpha=\alpha^ae_a$ in a basis and $g_{ab}=(e_a,e_b)$  (summation understood). If $\delta\alpha=0$ then
\[ \square\alpha=\alpha^a\delta\extd e_a+\del^a\alpha^b(\delta(e_ae_b)+e_a\vdash\extd e_b)+\del^a\del^b\alpha^c(e_a\vdash(e_be_c))\]
\end{corollary}
\proof The formula for $\delta(f\omega)$ follows immediately from the derivation property of $\extd$ and the first interior product property in  (\ref{inthodge}). The formula for $\square(f\omega)$ then follows from this and the Leibniz rule for $\extd$ as in \cite{Ma:rq}. These results then give the explicit formulae for $\alpha=\alpha^ae_a$. \endproof

Note concerning $\square\alpha$ that $(e_a\vdash(e_be_c))+ g_{bc}e_a-g_{ab}e_c=(e_ae_b)\dashv e_c$ in the classical case, which is  antisymmetric in $a,b$, while $L_\delta(\extd f,\omega)+\CL_{\extd f}\omega=2\nabla_{\extd f}\omega$ in the classical case as shown in \cite{Ma:rq}. Here $\nabla$ is the classical Levi-Civita connection referred back to a derivative along 1-forms via the metric.  The special case shown in Corollary~\ref{lbformula} is relevant to `Maxwell theory' where $F=\extd \alpha$ and Maxwell's equation $\delta F=J$ has a degree of freedom to change $\alpha$ by an exact form, which freedom can be reduced by fixing $\delta\alpha=0$. Maxwell's equation then becomes $\square\alpha=J$ where $J$ is required to be a  coexact `source'.

\subsection{Finite group case}

To give a concrete example we recall that any ad-stable subset of a finite group not containing the group identity defines an ad-stable ideal in $k(G)^+$ and hence a bicovariant calculus. There is a canonical basis of 1-forms $\{e_a\}$ labelled by the subset and relations $e_a f= R_a(f)e_a$ where $R_a(f)=f((\ )a)$ is right translation. This corresponds to crossed module action $e_a\ra f=f(a)e_a$. There is a natural choice of bi-invariant metric namely $g=\sum_a e_a\tens e_{a^{-1}}$ provided our subset is closed under inversion. Here the left coaction is trivial on $\Lambda^1$ and the right coaction is the crossed module one, namely $\Delta_Re_a=\sum_{g\in G}e_{gag^{-1}}\tens \delta_g$. The element $\theta=\sum_a e_a$ is similarly bi-invariant and makes the calculus inner (so $\extd=[\theta,\ \}$ as a graded commutator). The above results therefore give us a Hodge star on any finite group with bicovariant calculus stable under inversion and for which $B_-(\Lambda^1)$ has (up to scale) a unique top form which is bi-invariant. The braiding is $\Psi(e_a\tens e_b)=e_{aba^{-1}}\tens e_a$, from which the braided factorials can be computed. The metric sends $g(f^a)=e_{a^{-1}}$ where $\{f^a\}$ is a dual basis. The interior products given by the braided-derivatives
are 
\[ e_a\vdash(e_be_c)=g_{ab}e_c-g_{a,bcb^{-1}}e_b,\quad (e_ae_b)\dashv e_c=e_a g_{bc}-e_{aba^{-1}}g_{ac}\]
on evaluating against the super-braided coproduct via the metric. One also has $\extd\theta=0$ and in nice cases (when the associated quandle is locally skew and the field is characteristic zero) it was shown in \cite{MaRie} that $H^1(G)=k\theta$. If one similarly has $H^{n-1}(G)=k\theta^\sharp$ as an expression of Poincar\'e duality (which is often the case, including in the following example) then a source $J$ is coexact if and only if 
\[ \delta J=0,\quad \int_G (J,\theta)=0\]
by elementary arguments using (\ref{inthodge}), where $\int_G$ means a sum over the group. Here $\int:\Omega^n\to k$ is defined as the tensor product of $\int_G$ and $\int$ on $B_-(\Lambda^1)$ and $\int\extd \beta=0$ on any $n-1$-form $\beta$ under our assumption of a unique bi-invariant  top form ${\rm Vol}$ of degree $n$.

\begin{example}\label{FouS3} The standard 3D calculus on the permuation group $S_3$ on 3 elements is given by the conjugacy class of 2-cycles. We recall that $\Lambda^1=k-{\rm span}\{e_u,e_v,e_w\}$ is a $k(S_3)$-crossed module as above, where $a=u,v,w$ are the 2-cycles $u=(12), v=(23), w=uvu=vuv=(13)$. The minimal exterior algebra in this case is known to be a super version of the Fomin-Kirillov algebra\cite{FK, Ma:perm} with relations and exterior derivative
\[ e_u e_v + e_v e_w + e_w e_u=0,\quad e_v e_u + e_w e_v+ e_u e_w=0, \quad e_u^2=0,\quad \extd e_u+e_v e_w+ e_w e_v=0\]
and the two cylic rotations of these where $u\to v\to w\to u$. The dimensions in the different degrees are $\dim(\Lambda)=1:3:4:3:1$ so there is a unique top form up to scale, which we take as ${\rm Vol}=e_ue_ve_ue_w$. This is clearly central and one can check that it is also bi-invariant. This can be done noting that ${\rm Vol}=e_ue_ve_u\theta$ and computing
\[ \Delta_R(e_ue_ve_u)=e_ue_ve_u\tens(\delta_e+\delta_w)+e_we_ue_w\tens(\delta_u+\delta_{vu})+ e_ve_we_v\tens(\delta_{uv}+\delta_v).\] Hence we have a canonical Hodge star. 

The coevaluation element $\exp$ is a computation from Proposition~\ref{coev} which in a basis reads
\[ \exp=\sum_{m=0}^n \cdot[m,-\Psi]!^{-1}(e_{a_1}\tens \cdots\tens e_{a_m})\tens f^{a_m}\cdots f^{a_1}\]
where $\{e_a\}$ is our basis of $\Lambda^1$ and $\{f^a\}$ is a dual one and we sum repeated indices. This is in general but in our case and using the metric identification comes out as 
\[ \exp=1\tens 1+ \sum_a e_a\tens e_a+ e_v e_w\tens e_w e_v+ e_u e_w\tens e_w e_u - e_u e_v\tens e_u e_w-e_v e_u\tens e_v e_w\]
\[ \quad\quad +e_ue_ve_w\tens e_u e_v e_w+e_ve_we_u\tens e_v e_w e_u+ e_we_ue_v\tens e_we_ue_v-{\rm Vol}\tens{\rm Vol}\]
as one may check by verfying that this is a sum of basis and dual basis of each degree of $\Lambda^m\tens\Lambda^m$ paired via the metric. For example, 
\[ \<e_ue_v,e_we_v\>=\ev(e_u\tens e_v,e_w\tens e_v-e_u\tens e_w)=0\]
\[ \<e_ue_v,e_ve_w\>=\ev(e_u\tens e_v,e_v\tens e_w-e_u\tens e_v)=0\]
\[ \<e_ue_v,e_we_u\>=\ev(e_u\tens e_v,e_w\tens e_u-e_v\tens e_w)=0\]
\[ \<e_ue_v,e_ue_w\>=\ev(e_u\tens e_v,e_u\tens e_w-e_v\tens e_u)=-1\]
where we should remember that the $\ev$ pairing is nested starting on the inside. 

The resulting Hodge operator is then computed as
\[ \sharp 1=-{\rm Vol},\quad \sharp e_u=e_w e_u e_v,\quad \sharp e_v=e_u e_v e_w,\quad \sharp e_w=e_v e_w e_u\]
\[ \sharp (e_w e_u e_v)=-e_u,\quad \sharp(e_u e_v e_w)=-e_v,\quad \sharp(e_v e_w e_u)=-e_w,\quad \sharp {\rm Vol}=1\]
\[ \sharp(e_ue_v)=e_w e_u;\quad \sharp(e_a e_b)=e_{aba}e_a=\cdot\Psi(e_a\tens e_b)\]
We see that on the different degrees,
\[ \sharp^2|_{0,1,3,4}=-\id,\quad \sharp^3|_2=\id\]
so that $\sharp$ has order 6. The first of these illustrates Proposition~\ref{outer} while we note that $\sharp$ on degree 2 coincides with minus the braided-antipode $S$ of $B_-(\Lambda^2)$ (because this is braided-multiplicative along with an extra sign for the super case, and $S|_1=-\id$).  The cohomology for this calculus in characteristic zero is known to be $H^0=k, H^1=k, H^2=0, H^3=k, H^4=k$ and one can see that $\sharp$ is an isomorphism $H^m\isom H^{4-m}$ as expected, the cohomologies being spanned by $1, \theta=e_u+e_v+e_w$ in degrees 0,1 and their $\sharp$ in degrees 3,4. 

Next we compute $S$ on degree 3 as $S(e_ue_ve_w)=-\cdot\Psi(Se_u\tens S(e_ve_w))=\cdot\Psi(e_u\tens e_ue_v)=e_ue_we_u=-e_ue_ve_w$. A similar computation gives $S{\rm Vol}={\rm Vol}$ so that $S=(-1)^{D}$ on degree $D\ne 2$ as we saw in Proposition~\ref{outer}, while we have already observed that $S=-\sharp$ on degree 2.

We now compute $\sharp^\star:=\CF^\star\circ g$ as a check of our theory. The inverse braiding is $\Psi^{-1}(e_a\tens e_b)=e_b\tens e_{b^{-1}ab}$ in general and extends similarly to products, with the metric identification, we have
\[ \Psi_{sup}^{-1}\exp=1\tens 1- g+e_w e_v\tens e_w e_u+ e_w e_u\tens e_w e_v- e_u e_w\tens e_v e_w-e_v e_w\tens e_e e_w\]
\[ \quad -e_ue_ve_w\tens e_u e_v e_w- e_v e_w e_u\tens e_v e_w e_u- e_w e_u e_v\tens e_w e_u e_v-{\rm Vol}\tens{\rm Vol}.\]
Note that $e_ue_ve_w=e_we_ve_u$ and so forth using the relations. Integrating against gives $\sharp^\star=\sharp$ on degrees 0,1,3,4 as in Proposition~\ref{outer} while $\sharp^\star=\id$ on degree 2. The latter agrees with $\sharp^\star=\mu  (-1)^D S\sharp^{-1}$ in the proof of Proposition~\ref{outer}. 

Turning to applications we compute
\[ S\sharp 1=-{\rm Vol}, \quad S\sharp e_u=-e_we_ue_v,\quad S\sharp(e_ae_b)=-e_be_{bab},\quad S\sharp(e_we_ue_v)=e_u,\quad S\sharp{\rm Vol}=1\]
For the Laplace operator we note that  $\delta e_a=0$, for example 
\[ S\sharp\delta e_u= -\extd (e_w e_u e_v)=-\theta e_w e_u e_v-e_w e_u e_v \theta=-e_ue_we_ue_v- e_w e_u e_v e_u=0\]
using the relations. Hence by Corollary~\ref{lbformula} 
\[ \square|_0=\sum\del^a{}^2=-2\sum_a\del^a\]
which is (-2 times) the standard graph Laplacian for the corresponding Cayley graph on $S_3$. It is fully diagonalised as usual by the matrix elements of irreducible representations (the eigenvalues are 0, 6, 12 with eigenspaces of dimensions 1,4,1 respectively). We also have $\delta\alpha=g_{ab}\del^a\alpha^b$ and if this vanishes then 
\[ \square\alpha=(\sharp^{-1}\extd)^2\alpha=\sum (({1\over 2}\square+3)\alpha^a)e_a-(\alpha,\theta)\theta-\sum\del^a(\del^b+2)\alpha^{aba}e_b.\]  We note that the last term here only has contributions from $a\ne b$. The above expression is a short computation from Corollary~\ref{lbformula} using  \[ \delta(e_ae_b)=  e_b-e_{aba},\quad e_a\vdash(e_be_c)=\delta_{a,b}e_c-\delta_{aba,c}e_b\]
from which we see that
\[ \delta\extd e_a=3e_a-\theta,\quad \delta(e_ae_b)+e_a\vdash \extd e_b=e_b-2e_{aba}+\delta_{a,b}\theta.\]
Solving Maxwell theory in the form $\square\alpha=J$, the source $J$ has to be coexact. From the remarks above and $H^3(S_3)=k\theta^\sharp$, this is equivalent to 
\[ \delta J=0,\quad \int_{S_3} (J,\theta)=0\]
It is a useful check of our formula for $\square$ on $\Omega^1$ to see directly that when restricted to coexact forms its image indeed is again coexact. Moreover, by computer one finds  the same eigenspaces (each 4-dimensional) as in \cite{MaRai} with eigenvalues 3,6,9, so that up to an overall constant the Laplacian restricted to coexact 1-forms is the same in spite of the Hodge operators being rather different. Explicitly, 
\[ e_u-e_v, \quad e_u-\sum_a(\delta_{au}+\delta_{uau})e_a\]  
and their cyclic rotations under $u\to v\to w\to u$ have eigenvalue 3 (and along with their cyclic rotations add up to zero). Multiplying these by the sign function on $S_3$ gives eigenvectors of eigenvalue 9 while for the eigenvectors of eigenvalue 6 we can use the `point sources' in \cite{MaRai},
\[ J_x:=(3\delta_x-1)\theta+3\sum_a\delta_{xa}e_a,\quad x\in S_3\]
where three points that share a common node in the graph have a zero sum of their sources. These are related to the matrix elements $\rho_{ij}$ of the 2-dimensional representation. For example, if we work over $\R$ and $\rho(u)={\rm diag}(1,-1)$, $\rho(v)={1\over 2}{\rm diag}(-1,1)+{\sqrt{3}\over 2}\tau$ where $\tau$ is the transposition matrix, then 
\[ J_e=2\rho_{11}e_u+{\rm cyclic}, \quad J_{vu}-J_{uv}=2\sqrt{3}\rho_{21}e_u+{\rm cyclic}.\]
This solves the `Maxwell theory' on $S_3$ for this calculus by diagonalising $\square$ on coexact 1-forms.  
\end{example}

In \cite{MaRai}  and all other such models until now it has been  assumed that the Hodge operator should be designed to square to $\pm 1$, whereas our canonical Hodge operator in this example is order 6 and looks very different, but nevertheless gives the same reasonable Laplacians in degrees 0,1. Our construction also works for $S_4$ and $S_5$ with their 2-cycles calculus and can be analysed similarly, while higher $S_n$, $n>5$ are conjectured\cite{Ma:perm} to have infinite-dimensional $B_-(\Lambda^1)$. 

\section{Calculus and Hodge operator on coquasitriangular Hopf algebras}

Here we start with a new, braided-Lie algebra, approach to the construction of bicovariant $(\Omega^1,\extd)$ on quantum groups such as $k_q[G]$ for $G$ a complex semisimple Lie group. We recall that these are all {\em coquasitriangular} in that they come with a convolution-invertible map $\CR:A\tens A\to k$ obeying
\[ b\o a\o \CR(a\t\tens b\t)=\CR(a\o\tens b\o)a\t b\t \]
\[ \CR(ab\tens c)=\CR(a\tens c\o)\CR(b\tens c\t),\quad  \CR(a\tens bc)=\CR(a\o\tens c)\CR(a\t\tens b) \]
for all $a,b,c\in A$. This is just dual to Drinfeld's theory in \cite{Dri}, see \cite{Ma:bg, Ma:book}. We will need  the `quantum Killing form'
\[ \CQ(a\tens b)= \CR(b\o\tens a\o)\CR(a\t\tens b\t),\quad \forall a,b\in A\]
which obeys $\CQ(Sa\tens Sb)=\CQ(b\tens a)$ since $\CR$ is invariant under $S\tens S$.

The construction of differential calculi on a coquasitriangular Hopf algebra $A$ has its roots in R-matrix constructions from the 1990s but the following general construction builds on our recent treatment in \cite{MaTao}. It is shown there that $A$ is a left $A$-crossed module by
\begin{equation}\label{adLcrossed} \Ad_L(a)=a\o Sa\th\tens a\t,\quad a\la b=b\t \CR(b\o\tens a\o)\CR(a\t \tens b\th).\end{equation}
Then any subcoalgebra $\CL$ becomes a left $A$-crossed module by restriction and its dualisation $\CL^*$ in the finite-dimensional subcoalgebra case becomes a right $A$-crossed module. It is shown that the quantum Killing form regarded by evaluation on its first input as a map $\CQ:A^+\to \CL^*$ is a morphism of crossed modules. This gives:

\begin{proposition}cf \cite{MaTao} \label{Omega1coquasi} Let  $A$ be a coquasitriangular Hopf algebra and $\CL\subseteq A$ a nonzero finite-dimensional subcoalgebra. Then   $\Lambda^1={\rm image}(\CQ)$ and $\varpi=\CQ$  defines a bicovariant differential calculus  $\Omega^1$ on $A$. \end{proposition} 

In  \cite{MaTao} we used a version of this to naturally construct possibly non-surjective differential calculi with $\Lambda^1=\CL^*$, but we also see from this result that $\CL$ itself is the more fundamental object as starting point.

\subsection{Braided-Lie algebras} 

Our new approach is to start with a Hopf algebra $B$ in a braided category $\CC$ and find a `Lie algebra' for it. We then take its dual to define a calculus.

\begin{definition}\cite{Ma:blie} A left braided-Lie algebra is a coalgebra $\CL$  in a braided category,  together with a morphism $[\ ,\ ]:\CL\tens\CL\to \CL$ subject to the axioms shown diagrammatically in Figure~\ref{Lieax}. The  associated braiding $\tilde\Psi$ and braided-Killing form are also shown. \end{definition}

\begin{figure} \label{Lieax}
\[ \includegraphics[scale=.8]{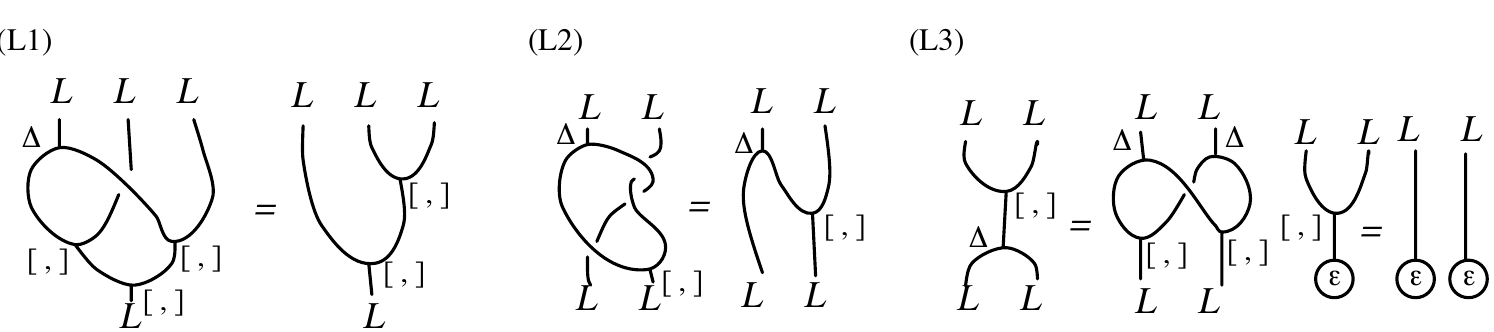}\]
\[   \includegraphics[scale=.5]{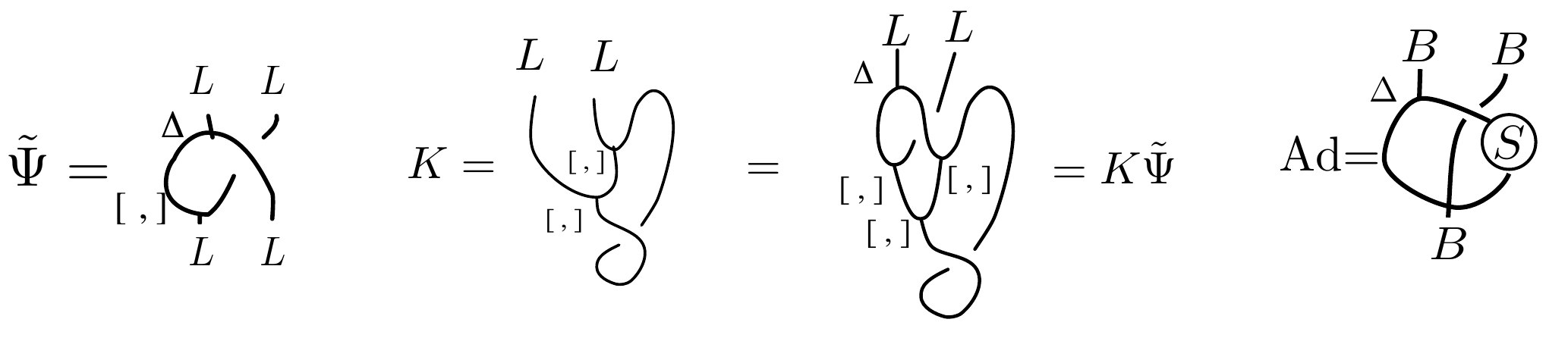} \]
\caption{Axioms of a braided-Lie algebra. Read down the page. We also recall the associated braiding $\tilde\Psi$ and braided-Killing form, and the adjoint action of a braided-Hopf algebra on itself.}
\end{figure}

A principal result in the case of an Abelian braided category is the construction of the braided-enveloping algebra $U(\CL)$ as a bialgebra. This is defined by the relations of commutativity with respect to the associated braiding $\tilde\Psi$. In the category of sets a braided-Lie algebra reduces to a quandle and this was used recently to prove the cohomology theorem for finite group bicovariant calculi\cite{MaRie}. Nondegeneracy of the Killing form also turns out to be an interesting characteristic related at one extreme to the Roth property of a finite group\cite{LMR}. The axioms themselves, however, were inspired by the properties of the braided adjoint action of a braided-Hopf algebra on itself as also recalled in Figure~3.

\begin{lemma}\cite{Ma:blie} If $B$ is a braided-Hopf algebra then $[\ ,\ ]=\Ad:B\tens B\to B$ the braided adjoint action obeys axiom (L1). If $B$ is cocommutative with respect to the braided-adjoint action  in the sense of \cite{Ma:bhop} (we say $B$ is {\rm Ad}-{\em cocommutative}) then (L2), (L3) are also obeyed. 
\end{lemma}

The first part was done in \cite{Ma:blie}. Braided cocommutativity with respect to a $B$-module is just the axiom (L2) when specialised to $[\ ,\ ]=\Ad$ and the proof that the adjoint action then obeys (L3) appeared in \cite[Prop A.2]{Ma:blin} in dual form (turn the diagrams there up-side-down). Clearly:

\begin{corollary} Any subcoalgebra $\CL\subseteq B$ of an Ad-cocommutative braided-Hopf algebra closed under $\Ad$ is a braided-Lie algebra by restriction. \end{corollary}

We next recall that if $A$ is coquasitriangular then there is a braided Hopf algebra version $B(A)$ of $A$ called its {\em transmutation}. This is also denoted $\und A$ and has the same coalgebra as $A$ but a modified product\cite{Ma:bg}
\begin{equation}\label{BAprod} a\bullet b=a\t b\th \CR(a\th \tens S b\o)\CR(a\o\tens b\t)=a\rz b\t \CR(a\co\tens S b\o)\end{equation}
 and lives in $\CM^A$ by $\Ad_R$. Its product is braided-commutative, 
\begin{equation}\label{BAcom} a\bullet b=b\th\bullet a\th\CR(Sb\t\tens a\o)\CR(b\fo\tens a\t)\CR(a\fo\tens b\fiv)\CR(a\fiv\tens Sb\o)\end{equation}
which can be written equivalently as
\begin{equation}\label{BAcomRE} \CR(b\o\tens a\o)a\t\bullet\CR(a\th\tens b\t) b\th=b\o\bullet \CR(b\t\tens a\o) a\t\CR(a\th\tens b\th)\end{equation}
while its braided-antipode is
\begin{equation}\label{BAS} \und S a=Sa\t\CR((S^2 a\th)Sa\o\tens a\fo)=(Sa\o)\rz\CR((Sa\o)\co\tens a\t).\end{equation}
Because we have a mix of both types of structure on the same vector space, we will be more careful to underline the braided versions where they are different.  

\begin{theorem}\label{Ablie} Let $A$ be coquasitriangular. Then $B(A)$ is $\und\Ad$-cocommutative and $[a,b]:=\und\Ad_a(b)=b\rz\CQ(a\tens b\co)$ for all $a,b\in A$ makes $A$ a braided-Lie algebra in the braided category $\CM^A$. By restriction, any subcoalgebra $\CL\subseteq A$ is a braided-Lie algebra in this category. In this case there is a surjection
\[ U(\CL)\twoheadrightarrow B(A)\]
of bialgebras in the category. 
\end{theorem}
\proof  We start by computing the left braided-adjoint action  by applying $\und S$ to $a\t$ and using the braiding to commute this past $b$ before multiplying up with respect to $\bullet$: 
\begin{eqnarray*} \und{\Ad}_a(b)&=& a\o\bullet b\rz\bullet(\und S a\t)\rz\CR((\und S a\t)\co\tens b\co)\\
&=& a\o\bullet \left( b\rz (\und S a\t)\rz\t\right)\CR((\und S a\t)\co\tens b\co\t)\CR(b\co\o\tens S(\und S a\t)\rz\o)\\
&=&a\o\rz b\rz\t (\und S a\t)\rz\th \CR((\und S a\t)\co\tens b\co\t)\CR(b\co\o\tens S(\und S a\t)\rz\o)\\
&&\quad\quad\quad\quad\CR(a\o\co\tens S(b\rz\o(\und S a\t)\rz\t)\\ 
&=&a\o\rz b\rz\t (Sa\t)\rz\th \CR((Sa\t)\co\t\tens a\th)\CR((Sa\t)\co\o\tens b\co\t)\\
&&\quad \CR(b\co\o\tens S(S a\t)\rz\o)\CR(a\o\co\tens S(b\rz\o(Sa\t)\rz\t))\\
&=&a\o\rz b\rz\t(Sa\t)\rz\t \CR((Sa\t)\co\tens a\th b\co\t)\\
&&\CR(a\o\co\t b\co\o\tens S(Sa\t)\rz\o)\CR(a\o\co\o\tens Sb\rz\o)\\
\end{eqnarray*}
where we use the definitions, the coaction properties and the multiplicativity property of $\CR$. We next unpack the adjoint coactions on $a$, and use multiplicativity of the last $\CR$ to give
\begin{eqnarray*} &=&a\th b\rz\th Sa\sev\CR(S^2 a\nine Sa\six\tens a\ten b\co\t)\CR(Sa\o a\fiv b\co\o\tens S^2 a\ei)\\
&&\CR(a\t\tens b\rz\t)\CR(a\fo\tens S b\rz\o)\\
&=& b\rz\t a\t Sa\sev\CR(S^2 a\nine Sa\six\tens a\ten b\co\t)\CR(Sa\o a\fiv b\co\o\tens S^2 a\ei)\\
&&\CR(a\th\tens b\rz\th)\CR(a\fo\tens S b\rz\o)\\
&=& b\rz\rz a\t Sa\six\CR(S^2 a\ei Sa\fiv\tens a\nine b\co\t)\CR(Sa\o a\fo b\co\o\tens S^2 a\sev)\\
&&\CR(a\th\tens b\rz\co)\\
&=& b\rz a\t Sa\six\CR(S^2 a\ei Sa\fiv\tens a\nine b\co\th)\CR(Sa\o a\fo b\co\t\tens S^2 a\sev)\\
&&\CR(a\th\tens b\co\o)\\
&=& b\rz a\t Sa\sev\CR(S^2 a\ten Sa\six\tens b\co\th)\nu^{-1}(a\ele)\CR(Sa\fiv\tens a\twe)\\
&&\CR(Sa\o a\fo\tens S^2 a\ei)\CR(b\co\t\tens S^2 a\nine)\CR(a\th\tens b\co\o)\\
&=& b\rz a\t Sa\sev\CR(S^2 a\ten Sa\six\tens b\co\th)\CR(Sa\fiv\tens S^2 a\ele)\nu^{-1}(a\twe)\\
&&\CR(Sa\o a\fo\tens S^2 a\ei)\CR(b\co\t\tens S^2 a\nine)\CR(a\th\tens b\co\o)\\
&=&\CR(a\o\tens Sa\ei) b\rz a\t Sa\sev\CR(Sa\fiv\tens S^2 a\twe)\nu^{-1}(a\thir)\\
&&\CR(S^2 a\ele Sa\six\tens b\co\th)\CR(b\co\t\tens S^2 a\ten)\CR(a\fo\tens S^2a\nine)\CR(a\th\tens b\co\o)\\
&=&\CR(a\o\tens Sa\ei) b\rz a\t Sa\sev\CR(Sa\six\tens S^2 a\ele)\nu^{-1}(a\thir)\\
&&\CR( Sa\fiv S^2 a\twe\tens b\co\th)\CR(b\co\o\tens S^2 a\nine)\CR(a\th\tens S^2a\ten)\CR(a\fo\tens b\co\t)\\
&=&\CR(a\o\tens Sa\six) b\rz a\t Sa\fiv \CR(a\th\tens S^2a\ei) \CR(Sa\fo\tens S^2 a\nine)\nu^{-1}(a\ele)\\
&&\CR( S^2 a\ten\tens b\co\t)\CR(b\co\o\tens S^2 a\sev)\CR(a\th\tens S^2a\ei)\\
&=&\CR(a\o\tens Sa\fo) b\rz a\t Sa\th\nu^{-1}(a\sev)\CR( S^2 a\six\tens b\co\t)\CR(b\co\o\tens S^2 a\fiv)\\
&=&\CR(a\o\tens Sa\t) b\rz \nu^{-1}(a\fiv)\CR( S^2 a\fo\tens b\co\t)\CR(b\co\o\tens S^2 a\th)\\
&=& b\rz \CR(a\t\tens b\co\t)\CR(b\co\o\tens a\o)=b\rz\CQ(a\tens b\co)=[a,b]
\end{eqnarray*}
where the 2nd equality is by quasicommutativity of $A$, the 3rd uses multiplicativity of $\CR$ to recognise $\Ad_R$ on $b\rz$. We then expand out by multiplicativity to recognise $\nu^{-1}(a)=\CR(S^2 a\o\tens a\t)$. This is known\cite{Ma:bg,Ma:book} to be convolution inverse to $\nu(a)=\CR(a\o\tens Sa\t)$ and to obey $\nu^{-1}(a\o)a\t=S^2 a\o\nu^{-1}(a\t)$, which we use to move to the right. The seventh equality uses multiplicativity of $\CR$ so that we can use quasicommutativity on $S^2a\ele Sa\six$ and the braid or Yang-Baxter equations on the last three factors to give the 8th equality. On this we use multiplicativity to cancel $a\fo Sa\fiv$ and obtain the 9th equality and two mutually inverse copies of $\CR$ for the 10th. We finally cancel $a\t Sa\th$ and move $\nu^{-1}$ to the left to cancel $\nu$. We then recognise the answer in terms of $\CQ$ and take this for our braided-Lie bracket. The Lemma tells us that we have (L1) for free. Next, we verify (L2) for $[\ ,\ ]=\und\Ad$ noting that (L2)  can be written in the form \[ \Psi\tilde\Psi=(\id\tens[\ ,\ ])(\Delta\tens\id)\]
where in our case $\Psi$ is the braiding of $\CM^A$ and 
\begin{eqnarray*}&&\nquad \tilde\Psi(a\tens b)= [a\o,b\rz]\tens a\t\rz\CR(a\t\co\tens b\co)\\
&&=b\rz\tens a\t\rz \CQ(a\o\tens b\co\o)\CR(a\t\co\tens b\co\t)\\
&&=b\t\tens a\fo \CR(((Sb\o)b\th)\o\tens a\o)\CR(a\t\tens ((Sb\o)b\th)\t)\\
&&\quad\quad\quad\quad\quad\CR((Sa\th)a\fiv\tens ((Sb\o)b\th)\th)\\
&&=b\t\tens a\t\CR(((Sb\o)b\th)\o\tens a\o)\CR(a\th\tens ((Sb\o)b\th)\t)\\
&&=b\rz\tens a\t \CR(b\co\o\tens a\o)\CR(a\th\tens b\co\t)
\end{eqnarray*}
using our result for $\und\Ad$. We compute
\begin{eqnarray*} 
&&\nquad \Psi\tilde\Psi(a\tens b)=a\t\rz\tens b\rz\rz\CR(b\rz\co\tens a\t\co)\CR(b\co\o\tens a\o)\CR(a\th\tens b\co\t)\\
&&=a\t\rz\tens b\rz\CR(b\co\o\tens a\t\co)\CR(b\co\t\tens a\o)\CR(a\th\tens b\co\th)\\
&&=a\t\t\tens b\rz\CR(b\co\o\tens a\o (Sa\t\o)a\t\th)\CR(a\th\tens b\co\t)\\
&&=a\o\tens b\rz\CR(b\co\o\tens a\t)\CR(a\th\tens b\co\t)=a\o\tens b\rz\CQ(a\t\tens b\co)=a\o\tens[a\t,b]
\end{eqnarray*}
as required, where we used the coaction properties of $\Ad_R$ and the multiplicativity property of $\CR$ to make a cancellation. The above Lemma then tells us that we get (L3) for free. These results then apply for an subcoalgebra $\CL\subseteq A$ since, due to the form of $\und\Ad$, we see that $\und\Ad(\CL\tens\CL)\subseteq \CL$, since $\Ad_R(\CL)\subseteq \CL\tens A$ because $\CL$ is a subcoalgebra (in other words a sub-coalgebra  of $A$ is also a subobject and hence a braided sub-coalgebra $\CL\subseteq B(A)$. For the last part, we can equivalently write 
\[ \tilde\Psi(a\tens b)=b\th\tens a\th \CR(Sb\t\tens a\o)\CR(b\fo\tens a\t)\CR(a\fo\tens b\fiv)\CR(a\fiv\tens Sb\o)\]
by expanding out our previous expression using the multiplicativity properties of $\CR$. Comparing with the braided-commutativity of $B(A)$ in (\ref{BAcom}) we see that $a\bullet b=\bullet\tilde\Psi(a\tens b)$ or the relations of $U(\CL)$. \endproof

\begin{proposition} \label{bliecoquasi}   When $\CL\subseteq A$ is finite-dimensional (1) $U(\CL)$ is Koszul dual to a right-handed quadratic version of the bicovariant calculus in Proposition~\ref{Omega1coquasi}. (2) The braided-Killing form is
\[ K(a,b)= \sum_i u(e_i\co\o)\CQ(a,e_i\co\t)\CQ(b,e_i\co\th) \<f^i,e_i\rz\>\]
where $\{e_i\}$ is a basis of $\CL$ and $\{f^i\}$ a dual basis and $u(a)=\CR(a\t\tens Sa\o)$.  \end{proposition}
\proof The braided-Killing form is obtained by reading down the diagram, as (summation understood)
\begin{eqnarray*}K(a,b)&=&\CQ(b,e_i\o)\CQ(a,e_i\rz\o)\<f^i\rz,e_i\rz\rz\rz\>\CR(e_i\rz\rz\co\tens f^i\co)\\
&=&\CQ(b,e_i\co)\CQ(a,e_i\rz\co)\CR(e_i\rz\rz\co\tens Se\rz\rz\rz\co)\<f^i,e_i\rz\rz\rz\rz\>\\
&=&\CQ(b,e_i\co\fo)\CQ(a,e_i\co\th)\CR(e_i\co\t\tens Se\co\o)\<f^i,e_i\rz\>
\end{eqnarray*}
as stated. For the remark about the dualisation we note that $A$ has a right crossed-module structure given by $\Ad_R$ and 
\[ a\ra b= a\t \CR(b\o\tens a\o)\CR(a\th\tens b\t),\quad \forall a,b\in A\]
and its crossed module braiding turns out, by similarly using the properties of $\CR$ as above, to coincide with the fundamental $\tilde\Psi$ for the braided-Lie algebra (as computed in the proof of Theorem~\ref{Ablie}). On the other hand this crossed module is the right handed version of (\ref{adLcrossed}) which dualized to give the crossed module structure defining the calculus in Proposition~\ref{Omega1coquasi}. This means that $U(\CL)$ is the Koszul or quadratic algebra dual of $\Lambda_{quad}$ (where we impose only the degree 2 relations of $B_-(\Lambda^1)$). The braided-Lie bracket and exterior derivative can also be related as part of a general theory of `quantum Lie algebras' in \cite{Wor} when $1\notin \CL$. Here every bicovariant calculus gives a quantum Lie algebra in the sense of \cite{Wor} and meanwhile (one can show that) every non-unital braided-Lie algebra $\CL$ gives a quantum Lie algebra by extending by $1$ and then taking the kernel of the counit. \endproof

The above theorem is a new result and is needed to complete the picture. In the special case where $\CL$ has a matrix coalgebra form on a basis $\{t^i{}_j\}$  (such data defines a matrix corepresentation of $A$) we recover the R-matrix braided-Lie algebra construction introduced in \cite{Ma:blie} but now as a corollary of the above.

\begin{corollary}cf\cite{Ma:blie}\label{matrixLie} Let $A$ be a coquasitriangular Hopf algebra and $t\in M_n(A)$ a matrix corepresentation. Then the matrix subcoalgebra $\CL=\{t^i{}_j\}$ has  braided-Lie bracket, categorical braiding and braided Killing form 
\[ [t^i{}_j,t^k{}_l]=t^{k_2}{}_{k_3}R^{-1}{}^{k_1}{}_{k_2}{}^i{}_{i_1}R^{k_3}{}_{k_4}{}^{i_1}{}_{i_2}R^{i_2}{}_{i_3}{}^{k_4}{}_l\tilde R^{i_3}{}_j{}^k{}_{k_1}\]
\[ \Psi(t^i{}_j\tens t^k{}_l)=t^{k_2}{}_{k_3}\tens t^{i_2}{}_{i_3} R{}^{i}{}_{i_1}{}^{k_1}{}_{k_2}R^{-1}{}^{i_1}{}_{i_2}{}^{k_3}{}_{k_4}R^{i_3}{}_{i_4}{}^{k_4}{}_l\tilde R^{i_4}{}_j{}^k{}_{k_1}\]
The braided enveloping algebra $U(\CL)$ is generated by the $\{t^i{}_j\}$ with new relations
\[ t^i{}_j\bullet t^k{}_l=t^{k_2}{}_{k_3}\bullet t^{i_2}{}_{i_3} R^{-1}{}^{k_1}{}_{k_2}{}^{i}{}_{i_1}R^{k_3}{}_{k_4}{}^{i_1}{}_{i_2}R^{i_3}{}_{i_4}{}^{k_4}{}_l\tilde R^{i_4}{}_j{}^k{}_{k_1}.\]
We sum over repeated indices in these expressions.
\end{corollary}
\proof We expand out $\CQ$ using the properties of $\CR$, then the above bracket can also be written explicitly as
\begin{equation}\label{bLiebracket} [a,b]=b\th \CR(Sb\t\tens a\o)\CR(b\fo\tens a\t)\CR(a\th\tens b\fiv)\CR(a\fo\tens Sb\o)\end{equation}
and the categorical braiding in $\CM^A$  is 
\begin{equation}\label{bLiebra} \Psi(a\tens b)= b\th\tens a\th\CR(a\o\tens b\t)\CR(Sa\t\tens b\fo)\CR(a\fo\tens b\fiv)\CR(a\fiv\tens Sb\o)\end{equation}
From these we immediately read off the expressions stated, where $R^i{}_j{}^k{}_l=\CR(t^i{}_j\tens t^k{}_l)$ and $\tilde R^i{}_j{}^k{}_l=\CR(t^i{}_j\tens St^k{}_l)$ is the `second inverse'. We likewise read off the relations of $U(\CL)$ from $\tilde\Psi$  or from (\ref{BAcom}) to give the result stated. In all cases we can move $\tilde R$ and another $R$ to the left hand side, for example the relations can be written compactly as  $R_{21}t_1\bullet R t_2=t_2\bullet R_{21}t_1R$ 
where the suffices refer to the position in $M_n\tens M_n$ with values in $U(\CL)$, also clear from (\ref{BAcomRE}). These are the relations of $B(R)$ \cite{Ma:exa,Ma:book}, the braided analog of the more familiar FRT bialgebra $A(R)$. \endproof

This derives the explicit R-matrix formulae needed to compute examples. This in turn recovers the 4D braided-Lie algebra of $k_q[SL_2]$ found in \cite{Ma:blie}:

\begin{example}\cite{Ma:blie}\label{LieSL2} For $A=k_q[SL_2]$ with $q^2\ne\pm 1$, its standard matrix coalgebra and rescaled generators
\[ t=\begin{pmatrix}\alpha &\beta\\ \gamma &\delta\end{pmatrix};\quad t=q^{-1}\alpha+q\delta,\quad z=\lambda^{-1}(\delta-\alpha),\quad x_+=\lambda^{-1}\beta,\quad x_-=\lambda^{-1}\gamma\]
where  $\lambda=1-q^{-2}$ (we use different symbols for the entries of $t^i{}_j$ to avoid confusion with the quantum group),  the nonzero braided-Lie brackets are
\[ [z,z]=q (2)_q\lambda z,\quad [t,t]=(2)_qt,\quad [t,\ ]=(q^3+q^{-3})\id,\quad [x_+,x_-]= z=-[x_-,x_+]\]
\[        [z,x_\pm]=\pm q^{\pm 1}(2)_q\,  x_\pm=-q^{\pm2}[x_\pm,z].   \]
Here $(2)_q=q+q^{-1}$ and we used Corollary~\ref{matrixLie} and the standard R-matrix for $SL_2$ with nonzero entries $R^1{}_2{}^2{}_1=q-q^{-1}, R^1{}_1{}^2{}_2=R^2{}_2{}^1{}_1=1, R^1{}_1{}^1{}_1=R^2{}_2{}^2{}_2=q$. The  braided Killing form is $[4,q^{-2}]/q^{10}$ times the nonzero values
\[ K(z,z)=(1+q^2),\quad K(t,t)=(1+q^6)(1+q^{-4}+\lambda),\quad K(x_+,x_-)=1=q^{-2}K(x_-,x_+). \]
The enveloping algebra $U(\CL)=B_q[M_2]$ is generated by $\alpha,\beta,\gamma,\delta$ with relations
\[ \beta\alpha=q^2\alpha\beta,\quad \gamma\alpha=q^{-2}\alpha\gamma,\quad \delta\alpha=\alpha\delta\]
\[ [\beta,\gamma]=\lambda\alpha(\delta-\alpha),\quad [\gamma,\delta]=\lambda\gamma\alpha,\quad [\delta,\beta]=
\lambda\alpha\beta\]
where $[\ ,\ ]$ at this point denotes commutator not Lie bracket. This is the algebra of $q$-deformed $2\times 2$ braided hermitian matrices which means that geometrically it should be thought of as $q$-Minkowski space\cite{Ma:exa,Ma:book}. There are two natural central elements, the braided determinant $\det_q=\alpha\delta-q^2\gamma\beta$ which should be thought of as the $q$-Lorentzian distance from the origin and $q$-trace ${\rm tr}_q=q^{-1}\alpha+q\delta=t$ which should be thought of as the `time' direction. In these variables (as opposed to the rescaled `Lie algebra' variables) the classical limit is commutative allowing us to think of this as a noncommutative geometry. Over $\C$ our braided-Lie algebra has a natural real form or $*$-involution $\alpha^*=\alpha, \beta^*=\gamma,\delta^*=\delta$ for real $q$, which fits with the mentioned geometric picture.  \end{example}

\subsection{Calculus and Hodge operator on $k_q[SL_2]$}

In the case of a coquasitriangular Hopf algebra $A$ with a generating matrix subscoalgebra $\{t^i{}_j\}$,  Proposition~\ref{Omega1coquasi} or dualization of Corollary~\ref{matrixLie} recovers a version of a known R-matrix construction\cite{Jurco} of quantum group covariant calculi. We let $\{E_\alpha{}^\beta\}$ be the standard basis of $M_n(k)$ and dual to the $\{t^i{}_j\}$ basis of $\CL$. This then becomes a right $A$-crossed module with
\begin{equation}\label{crossmat} \Delta_R E_\alpha{}^\beta=E_m{}^n\tens t^m{}_\alpha S t^\beta{}_n,\quad E_\alpha{}^\beta\ra t^a{}_b=E_m{}^n R^m{}_\alpha{}^a{}_c R^c{}_b{}^\beta{}_n\end{equation}
and
\[ \Lambda^1=M_n(k),\quad \varpi(a)=\CQ(a\tens t^\alpha{}_\beta)E_\alpha{}^\beta\]
defines the possibly non-surjective bicovariant calculus, which is, however typically surjective for the standard quantum groups $k_q[G]$ with $q$ generic. The $E_\alpha{}^\beta$ have bimodule relations
\[ E_\alpha{}^\beta t^a{}_b=t^a{}_c E_m{}^n R^m{}_\alpha{}^c{}_d R^d{}_b{}^\beta{}_n\]
and the above right covariance. The calculus has an inner form with \[ \extd t^a{}_b=t^a{}_c (R_{21}R)^c{}_b{}^\alpha{}_\beta E_\alpha{}^\beta-t^a{}_b\theta=[\theta,t^a{}_b],\quad \theta=E_\alpha{}^\alpha. \] 
The associated right crossed-module braiding on $\Lambda^1$ will be denoted $\tilde\Psi$  also (it is adjoint to the one for the braided-Lie algebra) and is computed from the right crossed module structure as
 \begin{eqnarray*}&&\nquad \tilde\Psi(E_\alpha{}^\beta\tens E_\gamma{}^\delta)=E{}_m{}^n\<E_\gamma{}^\delta, t^{j_1}{}_{j_2} \> \tens E_\alpha{}^\beta\ra t^m{}_{j_1} S t^{j_2}{}_n\\
 &&= E{}_m{}^n\tens E{}_p{}^q \<E_\gamma{}^\delta, t^{j_2}{}_{j_3} \> \<E_\alpha{}^\beta,t^{k_2}{}_{k_3}\>\CR(t^p{}_{k_1}\tens t^m{}_{j_1}S t^{j_4}{}_n)\CR(t^{j_1}{}_{j_2}S t^{j_3}{}_{j_4}\tens t^{k_3}{}_q)
 \end{eqnarray*}
 and expands out as 
\[ \tilde\Psi(E_\alpha{}^\beta\tens E_\gamma{}^\delta)=E_{j_2}{}^{j_3}\tens E_{k_2}{}^{k_3} \tilde R^{k_2}{}_{k_1}{}^{j_4}{}_{j_3}\, R^{k_1}{}_\alpha{}^{j_2}{}_{j_1}
R^{j_1}{}_\gamma{}^\beta{}_{k_4}R^{-1}{}^\delta{}_{j_4}{}^{k_4}{}_{k_3} \]
from which we see that $\tilde\Psi(E_\alpha{}^\beta\tens\theta)=\theta\tens E_\alpha{}^\beta$ so that, in particular, $\theta^2=0$ in $\Lambda_{min}=B_-(\Lambda^1)$. 

We now focus on $A=k_q[SL_2]$ where the smallest nontrivial irreducible is 2-dimensional, giving us $\Lambda^1=M_2(\C)$. We write basis $E_1{}^1=e_a, E_1{}^2=e_b, E_2{}^1=e_c, E_2{}^2=e_d$ and use the standard $SL_2$ R-matrix as in Example~\ref{LieSL2} to give the bimodule relations of the well-known 4D calculus first found in \cite{Wor}, 
\[ e_a
\begin{pmatrix}a&b\\ c&d\end{pmatrix}=\begin{pmatrix}qa&q^{-1} b\\
qc&q^{-1}d\end{pmatrix}e_a\]
\[  [e_b, \begin{pmatrix}a&b\cr c&d\end{pmatrix}]=q\lambda\begin{pmatrix}0&a\cr 0&c\end{pmatrix}e_a,
\quad [e_c, \begin{pmatrix}a&b\cr c&d\end{pmatrix}]=q\lambda\begin{pmatrix}b&0\cr d&0\end{pmatrix}e_a\]
\[ [e_d,\begin{pmatrix}a\cr c\end{pmatrix}]_{q^{-1}}=\lambda \begin{pmatrix}b\cr d\end{pmatrix}e_b,\quad 
[e_d,\begin{pmatrix}b\cr d\end{pmatrix}]_q=\lambda  \begin{pmatrix}a\cr c\end{pmatrix}e_c+q\lambda ^2 \begin{pmatrix}b\cr d\end{pmatrix}e_a,\]
where $[x,y]_q\equiv
xy-qyx$ and $\lambda=1-q^{-2}\ne 0$. The exterior differential is necessarily inner with $\theta=e_a+e_d$ which implies that
\[ \extd \begin{pmatrix}a\\ c\end{pmatrix}=\begin{pmatrix}a\\ c\end{pmatrix}((q-1)e_a+(q^{-1}-1)e_d)+\lambda \begin{pmatrix}b\\ d\end{pmatrix} e_b\]
\[ \extd \begin{pmatrix}b\\ d\end{pmatrix}=\begin{pmatrix}b\\ d\end{pmatrix}((q^{-1}-1+q\lambda^2)e_a+(q-1)e_d)+\lambda \begin{pmatrix}a\\ c\end{pmatrix} e_c.  \]
Note that we should scale $\extd$ or $\theta$ by $\lambda^{-1}$ in order to have the right classical limit but we have not done this in order to follow the general construction. The right coaction on  left-invariant 1-forms is 
\[ \Delta_R\theta=\theta\tens 1,\quad\Delta_R(-e_b,e_z,q^{-1}e_c)=(-e_b,e_z,q^{-1}e_c)\tens\begin{pmatrix} a^2 & (2)_qab & b^2\\
ca &1 + (2)_qbc & db \\ c^2 & (2)_qcd & d^2\end{pmatrix}\]
where $e_z:=q^{-2}e_a-e_d$ and the calculation is from $\Delta_R e_\alpha{}^\beta$ with the relevant R-matrix. We use the symmetric $q$-integers so that $(2)_q=q+q^{-1}$. The crossed module braiding comes out as 
\[ \tilde\Psi(e_a\tens \begin{pmatrix}e_a& e_b\cr e_c & e_d\end{pmatrix})= \begin{pmatrix}e_a& q^2 e_b\cr q^{-2}e_c & e_d\end{pmatrix}\tens e_a,\]
\[   \tilde\Psi(e_b\tens \begin{pmatrix}e_a& e_b\cr e_c & e_d\end{pmatrix})=\begin{pmatrix}e_a&  e_b\cr e_c & e_d\end{pmatrix}\tens e_b+ \lambda q^2 \begin{pmatrix}-e_b& 0 \cr  e_z  & e_b\end{pmatrix}\tens e_a\]
\[   \tilde\Psi(e_c\tens \begin{pmatrix}e_a& e_b\cr e_c &  e_d\end{pmatrix})=\begin{pmatrix}e_a&  e_b\cr e_c & e_d\end{pmatrix}\tens e_c+ \lambda \begin{pmatrix}e_c& - q^2 e_z  \cr  0 & - e_c\end{pmatrix}\tens e_a\]
\[   \tilde\Psi(e_d\tens e_a)=e_a\tens e_d+  \lambda^2 q^2 e_z\tens e_a-\lambda(e_b\tens e_c-e_c\tens e_b)\]
\[   \tilde\Psi(e_d\tens  e_b)=q^{-2}  e_b\tens e_d-  \lambda e_z\tens e_b  \]
\[ \tilde\Psi(e_d\tens e_c)=q^2 \lambda e_z\tens e_c+(q^4-1+q^{-2}) e_c \tens e_d+\lambda (q^4-1)e_c\tens e_z\]
\[ \tilde\Psi(e_d\tens e_d)=e_d\tens e_d +\lambda(e_b\tens e_c-e_c\tens e_b)-\lambda^2 q^2 e_z\tens e_a\]
from which one can see for example that $\tilde\Psi(e_i\tens\theta)=\theta\tens e_i$.  This then gives  the relations of $\Lambda=B_-(\Lambda^1)$ as usual Grassmann variables $e_a,e_b,e_c$ and\cite{Wor}
\[ e_a e_d+e_d e_a+\lambda e_c e_b=0,\quad e_d e_c+q^2e_c e_d+\lambda e_a e_c=0\]
\[ e_b e_d+q^2 e_d e_b+\lambda e_b e_a=0,\quad e_d^2=\lambda e_c e_b\]
or equivalently
\[ e_b e_z+q^2 e_z e_b=0,\quad  e_z e_c+q^2
e_c e_z=0\]
\[ 
e_z e_a+e_a e_z=\lambda e_c e_b,\quad e_z^2=(1-q^{-4})e_c e_b.\]
The exterior derivative is 
\[ \extd e_a=\lambda e_b e_c,\quad \extd e_c=\lambda q^2e_c
e_z,\quad \extd e_b=\lambda q^2e_z e_b,\quad \extd e_d= \lambda e_c e_b;\quad \extd
e_z=(1-q^{-4})e_b e_c.\]
As in degree 0, we note that $\lambda^{-1}\extd$ has the right classical limit not $\extd$ itself. The dimensions here in each degree are $\dim(\Lambda)=1:4:6:4:1$. The following is mostly known  e.g.\cite{Ma:ric} but we give a short proof as it is critical for us.

\begin{proposition}\label{4DSU2metric} For the above 4D calculus on $k_q[SL_2]$ with $q^2\ne\pm1$ there is a unique bi-invariant central metric
\[ g=e_c\tens e_b+ q^2 e_b\tens e_c + {q^3\over (2)_q}(e_z\tens e_z - \theta\tens\theta).\]
The inverse metric is
\[ (e_b,e_c)=1,\quad (e_c,e_b)=q^{-2},\quad (e_z,e_z)=q^{-3}(2)_q=-(\theta,\theta)\]
and the rest zero in this basis.  In the exterior algebra there is also up to scale a unique top form ${\rm Vol}=e_ae_be_ce_d=e_be_ce_ze_a$ and this is bi-invariant and central. Hence the super-braided Fourier transform applies and we have a Hodge star and interior products. 
\end{proposition}
\proof  Using the quantum group relations one has 
\[ \begin{pmatrix} a^2 & (2)_qab & b^2\\
ca &1 + (2)_qbc & db \\ c^2 & (2)_qcd & d^2\end{pmatrix}\begin{pmatrix} 0 & 0& -q^2\\ 0 & {q^3\over (2)_{q}}& 0\\
-1 & 0& 0\end{pmatrix} \begin{pmatrix} a^2 & (2)_qab & b^2\\
ca &1 + (2)_qbc & db \\ c^2 & (2)_qcd & d^2\end{pmatrix}^t=\begin{pmatrix} 0 & 0& -q^2\\ 0 & {q^3\over (2)_{q}}& 0\\
-1 & 0& 0\end{pmatrix}\]
which gives us the unique generically $q$-invariant element of the tensor square of the space spanned by $\{-e_b,e_z,q^{-1}e_c\}$ (the quadratic elements here generate  $k_q[SO_3]$). We can add to this  a multiple of $\theta\tens \theta$ since this is also invariant\cite{Ma:ric}, which we have now fixed so that $g$ is central. Thus for example,
\begin{eqnarray*} (e_c\tens e_b+ q^2 e_b\tens e_c)a & =& e_c\tens a e_b+ q^3\lambda e_b\tens be_a+q^2 e_b\tens a e_c\\
&=& a(e_c\tens e_b+ q^2 e_b\tens e_c)+q\lambda b e_a\tens e_b+ q^3\lambda be_b\tens e_a+q^4\lambda^2 a e_a\tens e_a\\
(e_z\tens e_z-\theta\tens\theta)a&=&-(1+q^{-2})(\lambda e_a\tens e_a+e_a\tens e_d+e_d\tens e_a)a\\
&=& a(e_z\tens e_z-\theta\tens\theta)-(1-q^{-4})(q^2\lambda a e_a\tens e_a+q b e_b\tens e_a+ q^{-1} b e_b\tens e_a)
\end{eqnarray*} using the above commutation relations. Comparing these we see that $[g,a]=0$. Similarly for the other generators of $k_q[SL_2]$. It is also clear that $\wedge(g)=0$. The inverse is immediate. For {\rm Vol} the element $e_be_ce_z$ is invariant again for reasons coming from the representation theory of $k_q[SO_3]$. As $\theta$ is also invariant, we know that $e_be_ce_z\theta$ is invariant and hence so is {\rm Vol} being a multiple of this. For centrality, we check for example
\[ ae_ae_be_ce_d=q^{-1}e_a a e_be_ce_d=e_ae_be_cq^{-1} ae_d=e_ae_be_ce_d a\]
discarding unwanted terms using the wedge product relations. 
\endproof

We now use $g$ to identify $\Lambda^1{}^*\isom \Lambda^1$ and compute $\sharp$ using $\int {\rm Vol}=1$.

\begin{proposition} For $k_q[SL_2]$ with its 4D calculus, $\mu=q^6$, $\sharp 1=q^6 {\rm Vol}$, 
\[ \sharp e_a=-q^4 e_ae_be_c,\quad \sharp e_b=-q^4 e_ae_be_d,\quad \sharp e_c=q^6 e_ae_ce_d,\quad \sharp e_d=q^4 e_be_ce_d+\lambda q^4 e_ae_be_c\]
\[ \sharp (e_ae_b)=-q^2 e_ae_b,\quad \sharp(e_ae_c)=q^4e_ae_c,\quad \sharp(e_ae_d)=q^2 e_be_c+\lambda q^4 e_a e_d\]
\[\sharp(e_be_c)=q^4 e_ae_d,\quad \sharp(e_be_d)=q^4 e_be_d+(1-q^4) e_ae_b,\quad \sharp(e_ce_d)=-q^2e_ce_d\]
\[ \sharp(e_ae_be_c)=-q^2 e_a,\quad\sharp(e_ae_be_d)=-q^2 e_b\quad \sharp(e_ae_ce_d)= e_c\quad \sharp(e_be_ce_d)=q^2e_d+\lambda q^2 e_a\]
and $\sharp {\rm Vol}=1$, where $\lambda=1-q^{-2}$ as above. Acting on degree $D$, this obeys 
\[ \sharp^2=q^6,\quad (D\ne 2);\quad (\sharp-q^4)(\sharp+q^2)=0,\quad (D=2).\]
\end{proposition}
\proof We first explicitly compute the exp element in the form
\[ \exp=1\tens 1+ g+ e_{i_1}e_{i_2}({}_2B){}^{-1}_{IJ}\tens e_{j_1} e_{j_2}+ e_{i_1}e_{i_2}e_{i_3} ({}_3B)^{-1}_{I J} \tens e_{j_1}e_{j_2}e_{j_3}+ e_1 e_2 e_3 e_4 ({}_4B)^{-1} \tens e_1 e_2 e_3 e_4\]
where $e_i, 1\le i\le 4$ refer in order to $e_a,e_b,e_c,e_d$  and $I=(i_1,j_2,\cdots,i_m)$ with $i_1<i_2\cdots<i_m$ labels of a basis of $\Lambda^m$ and 
\[ {}_mB{}_{IJ}=\<e_{i_1}\cdots e_{i_m},e_{j_1}\cdots e_{j_m}\>= \ev(e_{i_1}\tens e_{i_2}\cdots \tens e_{i_m}, [m, -\tilde\Psi]!(e_{j_1}\tens e_{j_2}\cdots e_{j_m}))\]
\[=g_{i_1 p_1}\cdots g_{i_mp_m}[m,-\tilde\Psi]!^{p_m\cdots p_2 p_1}_{j_1j_2\cdots j_m}\] 
In the last line refer operators to matrices, for example $[2,-\tilde\Psi](e_m\tens e_n)=e_p\tens e_q[2,-\tilde\Psi]^{pq}_{mn}$ and we remember the metric identification where $g_{ij}=(e_i,e_j)$. This is the general picture but with bases labelled in the classical way in the present example for generic $q$. We obtain
\[{}_1B^{-1}=\begin{pmatrix} -\lambda q^2 & 0 & 0 & -q^2 \\ 0 & 0 & q^2 & 0 \\ 0 & 1 & 0 & 0 \\ -q^2 & 0 & 0 & 0 \end{pmatrix}
    ,\quad  {}_2B^{-1}=q^2\begin{pmatrix} 0 & \lambda q^2 & 0 & 0 & 0 & -q^2 \\ \lambda & 0 & 0 & 0 & -1 & 0 \\ 0 & 0 & q^2 & 0 & 0 & 0 \\ 0 & 0 & 0 & 1 & 0 & 0 \\ 0 & -q^2 & 0 & 0 & 0 & 0 \\ -1 & 0 & 0 & 0 & 0 & 0 \end{pmatrix}\]
\[ {}_3B^{-1}=
q^4\begin{pmatrix}
 -\lambda  & 0 & 0 & -1 \\
 0 & 0 & q^2 & 0 \\
 0 & 1 & 0 & 0 \\
 -1 & 0 & 0 & 0 
\end{pmatrix}
,\quad {}_4B^{-1}=q^6\]
in the basis enumerations $12,13,14,23,24,34$ and $123,124,134,234$ respectively for the middle cases here. In particular, we see that $\mu=q^6$.  The matrix ${}_1B^{-1}$ here is the inverse of the matrix $g_{ij}$ in our basis and necessarily gives the coefficients of the metric $g\in \Lambda^1\tens\Lambda^1$, and we note also that $\mu=1/\det(g)$ in this basis. We then carefully integrate against this $\exp$, for example
\[ \sharp(e_be_d)=\int e_be_d e_ae_c \tens (\lambda q^2e_ae_b-q^2e_be_d)+\int e_be_d e_c e_d\tens (-q^2 e_ae_b)\]
where we read from the 2nd row of ${}_2B^{-1}$ for the terms in $\exp$ of the form $e_ae_c\tens\cdots$ and from the last row for terms of the form $e_ce_d\tens\cdots$. The other possibilities in our basis for the first tensor factor of $\exp$ have zero integral. We then evaluate the first displayed integral as $-q^2$ and the second integral as $\lambda$ on using the relations of the exterior algebra, to give $q^4e_be_d+(1-q^4)e_ae_b$ as stated. Integrating against $g$ is easier and gives $\sharp$ on degree 3.  We could now deduce $\sharp$ on degree 1  using Proposition~\ref{outer} that $\sharp^2=\mu$ on degrees $D\ne 2$ but one can also compute it similarly by integrating against ${}_3B^{-1}$ for a direct calculation and then verify $\sharp^2$. The polynomial identify for $\sharp$ on degree 2 is a direct calculation. \endproof

We see that $\sharp$ on degree 2 is not of finite order for generic $q$ but is a deformation of order 2. Indeed,  $\mu^{-1}=\det(g)$ (see the proof above) suggests a geometric normalisation to $\sharp'=\mu^{-{1\over 2}}\sharp=q^{-3}\sharp$ in our case,  then $\sharp'$ is involutive on degree $D\ne 2$ and obeys the standard $q$-Hecke relation 
\[ \sharp'{}^2=\id + (q-q^{-1})\sharp'\]
on degree 2. This is the same relation as obeyed by the braiding in the defining representation of the quantum group, which is also the braiding on the generators of the associated quantum plane. One can also compute $\sharp^\star$ directly and verify that it is given on degree 2 by $q^6S\sharp^{-1}$ as in Corollary~\ref{Fousup} (and otherwise equals $\sharp$). Here the braided antipode  on degree 2 is obtained from $S(e_ie_j)=-\cdot \tilde\Psi(e_i\tens e_j)$ as
\[ S(e_ae_b)=q^2e_ae_b,\quad S(e_ae_c)=q^{-2}e_ae_c,\quad S(e_ae_d)=e_ae_d-\lambda e_be_c\]
\[ S(e_be_c)=e_be_c +\lambda q^2 e_de_a,\quad S(e_be_d)=\lambda q^2 e_ae_b-e_de_b,\quad S(e_ce_d)=q^2e_ce_d.\]
One may similarly compute $S$ on degrees 3,4 
to  find $S=(-1)^D\id$ on all degrees $D\ne 2$ as must be the case by Proposition~\ref{outer}. Note that this feature of the antipode is not true for the outer degrees of all braided exterior algebras, see Example~\ref{Fouplane}. Finally, one can check that $[\sharp,S]=0$ on all degrees as it must by Corollary~\ref{Fousup}, which in degree 2 provides a very good cross-check of both the displayed $S$ and $\sharp$ computations. 

The cohomology for this calculus in characteristic zero and for generic $q$ is known to be $H^0=k, H^1=k, H^2=0, H^3=k, H^4=k$, which is the same as in Example~\ref{FouS3}. One can see that $\sharp$ is again an isomorphism $H^m\isom H^{4-m}$, the cohomologies being spanned by $1, \theta=e_a+e_d$ in degrees 0,1 and their $\sharp$ in degrees 3,4, where $\sharp\theta=-q^4e_be_c e_z$ deforms the classical 3-volume. Applications of this theory to $q$-electromagnetism will be considered elsewhere. Here we consider only the general result in Section~3.2 for the Laplacian on functions. We first recall from \cite{Ma:ric} that 
\begin{eqnarray}\label{extdSL2mon}
\extd(c^k b^n d^m)&=&  \left(q^{m+n-k}-1\right)\, c^k b^n d^m
e_d+
 \lambda q^{n}(k)_{q}\, c^{k-1} b^n\, d^{m+1}\, e_b \nonumber
\\
    &&\kern-20pt + \lambda q^{-k}\,
\left(q^{m-1}(m+n)_{q}\, c^{k+1} b^nd^{m-1}+(n)_{q}\,
c^kb^{n-1}d^{m-1}\right)e_c   \nonumber
\\
&&\kern-20pt +\lambda^2\, q \left((k+1)_{q}\, (m+n)_{q}\,  c^k b^n d^m
+
        q^{-m}(n)_{q}\, (k)_{q}\, c^{k-1}b^{n-1}d^m\right) e_a
        \nonumber
\\
    &&\kern-20pt +\left(q^{-m-n+k}-1\right)\, c^k b^n  d^m e_a 
\end{eqnarray}
where $(n)_q=(q^n-q^{-n})/(q-q^{-1})$.  Terms with negative powers of $c,b$ are treated as zero. 
This is only part of the algebra but there is a similar formula for the other part of the basis with $a$ in place of $d$ (in the classical limit this approach corresponds to patches where $d^{-1}$ and $a^{-1}$ respectively are adjoined). Writing 
\[ \extd f= (\del^bf)e_b+(\del^cf) e_c+(\del^zf) e_z+(\del^0 f)\theta\]
as the definition of our partial derivatives in this basis, it is shown in \cite{Ma:lap} that
 \[\del^0={q^2\lambda^2\over (2)_q}\Delta_q\ ,\] 
where \[  \Delta_q(c^k b^n d^m)=q^{-m}(k)_q(n)_q c^{k-1} b^{n-1} d^m+\left({k+n+m\over 2}\right)_q\left({k+n+m\over 2}+1\right)_q c^k b^n d^m  \]
is the naturally arising  $q$-deformed Laplace-Beltrami operator on $SL_2$.

\begin{proposition} For $q^2\ne\pm1$ and the quantum metric in Proposition~\ref{4DSU2metric} on the 4D calculus on $k_q[SL_2]$, the Hodge Laplacian on degree 0 is
\[ \square|_0=g_{ij}\del^i\del^j=2 q^{-1}\lambda^2\Delta_q\]
where $i,j$ are summed over our basis, for example $\{e_b,e_c,e_z,e_0\}$, and $g_{ij}=(e_i,e_j)$ are the metric coefficients.
\end{proposition}
\proof  Writing $e_a=(\theta+e_z)/(1+q^{-2})$ and $e_d=(\theta-q^2 e_z)/(1+q^2)$ we have 
\[ \del^z={1\over 1+q^{-2}}(\del^a-\del^d),\quad \del^0={1\over 1+q^{-2}}(\del^a+ q^{-2} \del^d),\]
in terms of the partial derivatives in our original basis read off from (\ref{extdSL2mon}). The former comes out as
out on $c^kb^nd^m$ as
\[ \del^z={\lambda\over (2)_q}(q^{m+n+2}(k)_q-q^{-k}(m+n)_q)+ {\lambda^2\over (2)_q}q^{2-m}(k)_q(n)_q S_c^- S_b^-\]
where $S_c^-$ lowers the degree of $c$ by 1, etc. A similar computation of $\del^a+q^{-2}\del^d$ gives $\del^0$ as stated. We next compute  on $c^k b^n d^m$ that
\begin{eqnarray*}&&\kern-20pt \del^b\del^c+q^{-2}\del^c\del^b = \del^b\lambda q^{-k+m-1}(m+n)_q S_c^+ S_d^-+ \del^b\lambda q^{-k}(n)_q S_b^- S_d^-+ q^{-2}\del^c\lambda q^n(k)_q S_c^- S_d^+\\
&&=2 \lambda^2 q^{n-k-1}(k)_q(n)_qS_c^-S_b^-+\lambda^2 q^{n-k-1+m}\left((k)_q(m+n+1)_q+(m+n)_q(k+1)_q\right)
\end{eqnarray*}
Next, we write for brevity
\[ A=q^{m+n+2}(k)_q-q^{-k}(m+n)_q,\quad B=\lambda q^2({k+n+m\over 2})_q({k+n+m\over 2}+1)_q\]
\[ T=\lambda  S_c^- S_b^- q^{2-m}(k)_q (n)_q\]
where $k,n,m$ are now the degree operators for the powers of $c,b,d$ respectively when acting on a monomial. Then $\del^z={\lambda\over(2)_q}(A+T)$ and $\del^0={\lambda\over(2)_q}(B+T)$ and 
\[ \del^z{}^2-\del^0{}^2={\lambda^2\over (2)_q^2}\left(A^2-B^2+2 T(A-B)\right)={\lambda\over (2)_q}(A+B+ 2T)\left(1-q^{m+n-k}\right)\]
on noting that $A-B={(2)_q\over\lambda}(1-q^{m+n-k})$ commutes with $T$ (since the latter changes both $k,n$ equally and does not change $m$). Putting in these results and the value of $A+B$ we obtain $\del^b\del^c+q^{-2}\del^c\del^b+q^{-3}(2)_q(\del^z{}^2-\del^0{}^2)= 2q^{-1}\lambda^2\Delta_q$. This can also be used to expresses $\del^0$ or $\Delta_q$ in terms of $\del^b,\del^c,\del^z$. 

It remains to check that $\delta e_i=0$ so that Corollary~\ref{lbformula} applies and this is the Hodge Laplacian. Indeed $\extd (S\sharp)^{-1}e_i=-q^{-6}\extd\sharp e_i$ and $\extd\sharp e_a=0$ since
\[ \theta\sharp e_a=-q^4\theta e_ae_be_c=-q^4e_de_ae_be_c=q^4e_a e_be_ce_d=q^4e_ae_be_c\theta=-(\sharp e_a)\theta\]
using the relations of the exterior algebra. Meanwhile one finds more strongly that
\[ \theta \sharp e_i=(\sharp e_i)\theta=0;\quad i=b,c,z\]
where $\sharp e_z=q^4 e_be_c\theta$, again using the relations.  \endproof

Note that we have been working algebraically but, over $\C$, our constructions are compatible with the  $*$-involution corresponding to the compact real form $SU_2$.


\begin{thebibliography}{999}
\bibitem{Baz}
Y. Bazlov,  Nichols-Woronowicz algebra model for Schubert calculus on Coxeter groups, J. Algebra 297 (2006) 372--399

\bibitem{BKLT} Y. Bespalov, T. Kerler, V. Lyubashenko and V. Turaev, Integrals for braided Hopf algebras, J. Pure Appl. Algebra 148 (2000), 113--164

\bibitem{Brz}
T. Brzezinski, Remarks on bicovariant differential calculi and exterior Hopf algebras, Lett. Math. Phys. 27 (1993) 287--300.

\bibitem{Cas}
L. Castellani, R. Catenacci and P.A. Grassi. The geometry of supermanifolds and new supersymmetric actions 
Nucl.Phys. B 899 (2015) 112-148.

\bibitem{Dri}
V.G. Drinfeld, Quantum Groups, in {\em Proc. of the ICM}, AMS (1987)

\bibitem{FK} S. Fomin and A.N. Kirillov, Quadratic algebras, Dunkl elements, and Schubert calculus, Adv.
Geom, Progr. Math. 172 (1989), 147--182

\bibitem{GomMa:ele}X. Gomez and S.Majid, Noncommutative cohomology and electromagnetism on $\C_q[SL2]$ at roots of unity, Lett. Math. Phys. 60 (2002) 221--237

\bibitem{GomMa} X. Gomez and S.Majid, Braided Lie algebras and bicovariant differential calculi over coquasitriangular Hopf algebras, J. Algebra 261 (2003) 334--388

\bibitem{Heck}I. Heckenberger, Hodge and Laplace-Beltrami operators for bicovariant differential calculi on quantum groups, Compositio Mathematica 123 (2000) 329--354

\bibitem{Jurco}  B. Jurco, Differential calculus on quantized simple Lie groups, Lett. Math. Phys. 22 (I991) 177--186

\bibitem{KemMa} A. Kempf and S. Majid, Algebraic q-integration and Fourier theory on quantum and braided spaces, J. Math. Phys. 35 (1994) 6802--6837

\bibitem{Kus} J. Kustermans, G. J. Murphy and L. Tuset, Quantum groups, differential calculi and the eigenvalues of the Laplacian, Trans. Amer. Math. Soc. 357 (2005) 4681--4717

\bibitem{LMR} J. Lopez Pena, S. Majid and K. Rietsch, Lie theory of finite simple groups and the Roth property, 39pp., arXiv:1003.5611 (math.QA)

\bibitem{Lus} G. Lusztig, {\em Introduction to Quantum Groups} (Birkhauser, 1993).

\bibitem{LyuMa} V. Lyubashenko and S. Majid, Braided groups and quantum Fourier transform, J. Algebra. 166 (1994) 506--528

\bibitem{Lyu} V. Lyubashenko, Modular transformations for tensor categories, J. Pure Appl. Algebra 98
(1995), 279--327

\bibitem{Ma:self} S. Majid, The self-representing Universe, in {\em Mathematical Structures of the Universe}, eds. M. Eckstein, M. Heller, S. Szybka, Copernicus Center Press (2014) 357 -- 387

\bibitem{Ma:book}
S. Majid, {\em Foundations of Quantum Group Theory}, C.U.P. (2000) 640pp

\bibitem{Ma:bhop} S. Majid, Algebras and Hopf algebras in braided categories, in Lec. Notes Pure and Applied Maths 158 (1994) 55--105. Marcel Dekker

\bibitem{Ma:exa} S. Majid, Examples of braided groups and braided matrices, J. Math. Phys. 32 (1991) 3246--3253

\bibitem{Ma:bg} S. Majid, Braided groups, J. Pure Applied Algebra 86 (1993) 187--221

\bibitem{Ma:cro} S. Majid, Cross products by braided groups and bosonization, J. Algebra 163 (1994) 165--190

\bibitem{Ma:skl}
S. Majid, Braided matrix structure of the Sklyanin algebra and of the quantum Lorentz group, Comm. Math. Phys. 156(1993) 607--638

\bibitem{Ma:fre} S. Majid, Free braided differential calculus, braided binomial theorem and the braided exponential map, J. Math. Phys. 34 (1993) 4843--4856

\bibitem{Ma:blin} S. Majid, Quantum and braided linear algebra, J. Math. Phys. 34 (1993) 1176--1196

\bibitem{Ma:blie} S. Majid, Quantum and braided Lie algebras, J. Geom. Phys. 13 (1994) 307--356

\bibitem{Ma:dbos}
S. Majid, Double bosonisation of braided groups and the construction of $U_q(g)$, Math. Proc. Camb. Phil. Soc.125 (1999) 151--192

\bibitem{Ma:dcalc}
S. Majid, Classification of differentials on quantum doubles and finite noncommutative geometry, Lect. Notes Pure Appl. Maths 239 (2004) 167--188, Marcel Dekker

\bibitem{Ma:perm} S. Majid, Noncommutative differentials and Yang-Mills on permutation groups $S_N$, Lect. Notes Pure Appl. Maths 239 (2004) 189--214, Marcel Dekker

\bibitem{Ma:ric} S.Majid,  Noncommutative Ricci curvature and Dirac operator on $\C_q[SL2]$ at roots of unity, Lett. Math. Phys. 63 (2003) 39--54

\bibitem{Ma:lap} S. Majid, Algebraic approach to quantum gravity III: noncommutative Riemannian geometry, in {\em Mathematical and Physical Aspects of Quantum Gravity}, eds. B. Fauser, J. Tolksdorf and E. Zeidler, Birkhauser (2006) 77--100

\bibitem{Ma:rq} S. Majid, Reconstruction and quantisation of Riemannian structures, 40pp. arXiv:1307.2778 (math.QA)

\bibitem{MaRai} S. Majid and E. Raineri, Electromagnetism and gauge theory on the permutation group $S_3$, J. Geom. Phys. 44 (2002) 129--155

\bibitem{MaRie} S. Majid and K. Rietsch, Lie theory and coverings of finite groups, J. Algebra, 389 (2013) 137--150

\bibitem{MaTao}
S. Majid and W.-Q. Tao, Duality for generalised differentials on quantum groups, J. Algebra, 439 (2015) 67--109 

\bibitem{Man}Yu.  Manin, {\em Quantum groups and noncommutative geometry}, Centre De Recherches Mathematiques (1988)

\bibitem{Nic} W. Nichols, Bialgebras of type I. Comm. Algebra 15 (1978) 1521--1552


\bibitem{Rad}
D. Radford, The structure of Hopf algebras with a projection, J. Algebra, 92 (1985), 322--347


\bibitem{Wor}
S. L. Woronowicz, Differential calculus on compact matrix pseudogroups (quantum groups), Comm. Math. Phys. 122 (1989), 125--170

\bibitem{Whi} J.H.C. Whitehead. Combinatorial homotopy, II. Bull. Amer. Math. Soc., 55:453--496, 1949.
51

\bibitem{Yet} D.N. Yetter, Quantum groups and representations of monoidal categories, Math. Proc. Camb. Phil. Soc., 108 (1990) 261--290

\end{thebibliography}
\end{document}